\documentclass[alpha-refs]{wiley-article}

\usepackage{lineno}
\usepackage{multirow}
\usepackage[usestackEOL]{stackengine}
\usepackage{amssymb}
\usepackage{pifont}
\newcommand{\cmark}{\ding{51}}%
\newcommand{\xmark}{\ding{55}}%

\papertype{Submitted to QJRMS, 2019}
\paperfield{}

\title{Numerical methods for entrainment and detrainment in the multi-fluid Euler equations for convection}

\abbrevs{ABC, a black cat; DEF, doesn't ever fret; GHI, goes home immediately.}

\author[1\authfn{1}]{William A McIntyre}
\author[1\authfn{1}]{Hilary Weller}
\author[1\authfn{2}]{Christopher E Holloway}

\contrib[\authfn{1}]{Equally contributing authors.}

\affil[1]{Department of Meteorology, University of Reading, UK}

\corraddress{Author One PhD, Department, Institution, City, State or Province, Postal Code, Country}
\corremail{h.weller@reading.ac.uk}

\fundinginfo{Funder One, Funder One Department, Grant/Award Number: 123456, 123457 and 123458; Funder Two, Funder Two Department, Grant/Award Number: 123459}

\begin{document}

\maketitle

\begin{abstract}
Convection schemes are a large source of error in global weather and climate models, and modern resolutions are often too fine to parameterise convection but are still too coarse to fully resolve it. Recently, numerical solutions of multi-fluid equations have been proposed for a more flexible and consistent treatment of sub-grid scale convection, including net mass transport by convection and non-equilibrium dynamics. The technique involves splitting the atmosphere into multiple fluids. For example, the atmosphere could be divided into buoyant updrafts and stable regions. The fluids interact through a common pressure, drag and mass transfers (entrainment and detrainment). Little is known about the numerical properties of mass transfer terms between the fluids. 
We derive mass transfer terms which relabel the fluids and derive numerical properties of the transfer schemes, including boundedness, momentum conservation and energy conservation. Numerical solutions of the multi-fluid Euler equations using a C-grid are presented using stable and unstable treatments of the transfers on a well-resolved two-fluid dry convection test case. We find two schemes which are conservative, stable and bounded for large timesteps, and maintain their numerical properties on staggered grids.

\end{abstract}

\section{Introduction}
The modelling of atmospheric convection is at the forefront of current meteorological research due to the large errors caused by convection schemes in atmospheric models \cite[eg][]{arakawa2004,lean2008,yano2004,yano2018}. We have reached the ``grey zone'' in which improved computational power allows convection to be partially resolved but cannot yet be explicitly simulated \cite[]{gerard2009}. Many convection schemes assume the scales of convection are small relative to the dynamical flow and that there is no net mass transport due to convection \cite[eg][]{arakawa1974}. This assumption is increasingly unrealistic at finer resolutions \cite[]{kwon2017}, since mass transport by convection could be larger than other mass fluxes once a convective cell is close to the grid-scale. This possibility is commonly ignored in convection schemes \cite[such as][]{arakawa1974,gregory1990,lappen2001} and so is non-equilibrium dynamics \cite[for example][]{kain1990}. Although many convection schemes incorporate some of these aspects \cite[such as][]{gerard2005,kuell2008}, few of them offer a consistent treatment of resolved and sub-grid convection. The conditional filtering (or conditional averaging) technique has been proposed for convection modelling due to the possibility of modelling convection over all resolutions, whilst also representing net mass transport by convection and non-equilibrium dynamics \cite[]{thuburn2017,weller2018}.  \\

Conditional filtering involves dividing space into various fluids. In the case of convection, one could label fluid 0 as regions of neutrally buoyant air, fluid 1 as convective updrafts and fluid 2 as downdraft regions. One then calculates a convolution with a space-dependent filter - which may be a volume average or a more complicated filter such as a Guassian. Each filtered fluid has its own properties and prognostic variables such as volume fraction, density, temperature and velocity \citep{thuburn2017}, and the equations of motion are solved for each fluid individually. As the scheme allows for the advection of any fluid to neighbouring cells, net mass transport by convection can occur in which the properties of the convective mass are transported based on the local dynamics. \\

Conditional filtering is used in other fields of science and engineering  \cite[]{dopazo1977,baer1986,mechitoua2003,guelfi2007}. \cite{lappen2001} used the technique to model cumulus convection, using an updraft fluid and a stable fluid. However, the stable fluid in that study was assumed to subside in the same column as the updraft, meaning the scheme does not incorporate net vertical mass transport by convection. More recently, \cite{thuburn2017} and \cite{tan2018} have described how to conditionally filter the fully compressible Euler equations with the aim of representing sub-grid scale convection. Subsequent studies have built upon these foundations including \cite{thuburn2018} who investigate the conservation properties and normal modes of the equations, \cite{thuburn2019} who use the method for a 2-fluid single-column convective boundary layer scheme, and \cite{weller2018} who formulate a numerical solution of the multi-fluid compressible Euler equations. Thus far, little is known about the numerical properties of solutions to the multi-fluid equations. \cite{stewart1984} and \cite{thuburn2019} note that the multi-fluid Euler equations are ill-posed when sub-filter terms are ignored. This property is confirmed by \cite{weller2018}, as drag or mixing between the fluids is necessary to prevent the fluid properties unphysically diverging from each other. However, the numerical properties of transfer terms between fluids has received little attention. \\

Transfer terms that exchange mass and other properties between fluids are crucial for formulating a parameterisation of convection using conditional filtering. These transfer terms will be equivalent to entrainment (including cloud base entrainment) and detrainment which may be adapted from existing frameworks such as \cite{arakawa1974}, \cite{betts1993}, \cite{neggers2002} and \cite{siebesma2007}. Fluid transfer terms are given in \cite{thuburn2018} and \cite{weller2018} in terms of transfer rates. \cite{weller2018} also propose a numerical scheme for the fluid transfers, but only one transfer scheme is considered in which the numerical treatment of the mass transfer is explicit (and the momentum/temperature transfer is treated implicitly) which may not be suitable for all transfer rates. This motivates us to present more mass transfer schemes and analyse their numerical properties to obtain the most desirable numerical solutions. \\

In this study, we analyse the numerical properties of the transfer terms between fluids for the multi-fluid compressible Euler equations (defined in section \ref{section_swe_single_fluid}). We formulate 20 possible numerical schemes and analyse their properties including conservation, boundedness and stability in section \ref{section_transfers}. We then apply the transfer terms to well-resolved two-fluid dry convection test cases in section \ref{section_rising_bubble}.

\section{Governing equations}
\label{section_swe_single_fluid}
The multi-fluid compressible Euler equations are derived in \cite{thuburn2017} and we will be using the notation convection from \cite{weller2018}. We have three equations for each fluid including the continuity equation,

\begin{equation}
\label{continuity_equation}
\frac{\partial \eta_i}{\partial t} + \bm{\nabla} . (\eta_i \bm{u}_i) = \underbrace{ \sum_{j \neq i} (\eta_j S_{ji} - \eta_i S_{ij}) }_\text{\clap{Mass transfers~}},
\end{equation}

the potential temperature equation,

\begin{equation}
\label{temperature_equation}
\frac{\partial \theta_i}{\partial t} + \bm{u}_i . \bm{\nabla} \theta_i = \underbrace{ \sum_{j \neq i} \left( \frac{\eta_j}{\eta_i} S_{ji} (\theta_j - \theta_i) \right) }_\text{\clap{Transfer mean temp.~}} - \underbrace{ \sum_{j \neq i} H_{ij} }_\text{\clap{Heat transfer~}},
\end{equation}

and momentum equation,

\begin{equation}
\label{momentum_equation}
\frac{\partial \bm{u}_i}{\partial t} + \bm{u}_i . \bm{\nabla} \bm{u}_i = \bm{g} - \underbrace{c_p \theta_i \bm{\nabla} \pi }_\text{\clap{Pressure gradient~}} \ + \ \underbrace{ \sum_{j \neq i} \left( \frac{\eta_j}{\eta_i} S_{ji} (\bm{u}_j - \bm{u}_i) \right) }_\text{\clap{Transfer mean velocity.~}} - \underbrace{ \sum_{j \neq i} \bm{D}_{ij} }_\text{\clap{Drag~}},
\end{equation}

where $i$ is the label for fluid $i$. Our prognostic variables are the fluid mass per unit volume ($\eta_i$), the fluid potential temperature ($\theta_i$) and the fluid velocity ($\bm{u}_i$). We have defined $\eta_i \equiv \sigma_i \rho_i$, where $\sigma_i$ is the fluid volume fraction and $\rho_i$ is the density of fluid $i$. $S_{ij}$ is the unidirectional mass transfer rate from fluid $i$ to fluid $j$ ($S_{ij} \geq 0$). $c_p$ is the heat capacity of dry air at constant pressure and $\bm{g}$ is the gravitational acceleration. $\pi \equiv p^\kappa/p_0^\kappa$ is the Exner pressure where $p$ is the pressure, $p_0$ is a reference pressure, $\kappa=R/c_p$ and $R$ is the gas constant of dry air. $H_{ij}$ is the heat transfer between fluids $i$ and $j$ (which does not exchange mass between fluids) and $\bm{D}_{ij}$ is the drag between fluids $i$ and $j$ - these exchange terms will not be used in this study. Additionally, the equation of state for dry air is used to relate the pressure and fluid temperatures:

\begin{equation}
p_0 \pi^{\frac{1-\kappa}{\kappa}} = R \sum_i \eta_i \theta_i.
\end{equation}

The total energy of the multi-fluid system is given by 

\begin{equation}
E = E_P + E_I + E_K
\end{equation}

where $E_P$ is the potential energy, $E_I$ is the internal energy and $E_K$ is the kinetic energy, defined respectively as:

\begin{align}
E_P &= \sum_i \eta_i |\bm{g}| z,\\
E_I &= \sum_i \eta_i \theta_i c_v \pi,\\
E_K &= \sum_i \frac{1}{2} \eta_i \bm{u}_i . \bm{u}_i,
\end{align}

where $z$ is the height coordinate and $c_v=\frac{c_p}{\gamma}$ is the heat capacity of dry air at constant volume and $\gamma$ is the heat capacity ratio. These energies will be used to assess the numerical stability of the transfer schemes.

\section{Fluid transfer schemes}
\label{section_transfers}
In convection modelling, entrainment and detrainment are both mass exchanges between the updraft and surrounding environment. It is therefore important that the numerical implementation of these transfer terms for a multi-fluid system have accurate conservation properties and that they do not produce new extrema. For accuracy, mass and momentum should be conserved. For stability, the fluid mass ($\eta_i$) must remain positive and velocities should be bounded. For accuracy and stability, potential and internal energy should also be conserved and the kinetic energy should not increase (resolved kinetic energy decreases when two fluids of differing velocity mix). In this section, we demonstrate the conservation and boundedness properties of the mass transfer terms for the multi-fluid equations and present solutions of the multi-fluid Euler equations with these transfers. 

\subsection{Notation and numerics}
In our governing equations, we have assumed that mass transferred between fluids will take its associated mean properties from the original fluid. Refinement of the transfer terms to incorporate sub-filter-scale variation will not form part of this study. For a mean fluid property $\phi_i \in [\theta_i, \bm{u}_i]$, the governing equations can be generalised as

\begin{equation}
\label{momentum_partitioned_advective_transfers}
\frac{\partial \phi_i}{\partial t} + \bm{u}_i.\bm{\nabla} \phi_i = F_i + \sum_{j \neq i} S_{ji} \frac{\eta_j}{\eta_i} \left[ \phi_j - \phi_i \right] ,
\end{equation}

where $F_i$ contains right-hand-side terms such as the pressure gradient term. Applying this to the temperature and momentum equations we get:

\begin{itemize}
	\item Momentum equation: $\bm{\phi}_i=\bm{u}_i$, $\bm{F}_{\bm{u}i} = - c_p \theta_i \bm{\nabla} \pi + \bm{g}$.
	\item Temperature equation: $\phi_i= \theta_i$, $\bm{F}_{\theta i} = 0$. 
\end{itemize}

We will assume that the transfer terms are operator split such that other processes (advection and $F_i$) act on the prognostic variables first, followed by the transfers: 

\begin{align}
\phi_i^m &= \phi_i^n - (1-\alpha) \Delta t \left[ \bm{u}_i.\bm{\nabla} \phi_i - F_i \right]^{n} + \alpha \Delta t \left[ \bm{u}_i.\bm{\nabla} \phi_i - F_i \right]^{m}, \\
\phi_i^{n+1} &= \phi_i^m + \Delta t \sum_{j \neq i} S_{ji} \frac{\eta_j}{\eta_i} \left[ \phi_j - \phi_i \right],
\end{align}

where $n$ is the time-level ($t = n \Delta t$), $m$ is the time-level after applying the advection and $F_i$ terms and $\alpha$ is the Crank-Nicolson off-centering coefficient. The mass transfers are then based on the most up-to-date states ($m$) rather than the previous time-level ($n$). This allows the transfer terms to be independent of the numerical properties of the advection and $F_i$ terms. For each momentum equation, the total momentum should be equal before ($m$) and after ($n+1$) mass transfer such that: $\sum_i \eta_i^{n+1} \bm{u}_i^{n+1} = \sum_i \eta_i^{m} \bm{u}_i^{m}$. The internal energy should also be conserved by the temperature equation transfers: $c_v \pi \sum_i \eta_i^{n+1} \theta_i^{n+1} = c_v \pi \sum_i \eta_i^{m} \theta_i^{m}$.

\subsection{Transferring mass}
\label{section_transfer_continuity}
When transferring mass, we must ensure that mass is conserved and all $\eta_i$ remain positive.  The mass in each fluid at the end of the timestep ($\eta_i^{n+1}$) is given by the mass after advection ($\eta_i^{m}$) plus the discretised transfer term integrated over time $\Delta t$:

\begin{equation}
\label{mass_transfer_formulation}
\begin{split}
\eta_0^{n+1} &= \eta_0^{m} - \Delta t \left[ (1-\alpha_C)\ \eta_0^{m} + \alpha_C\ \eta_0^{n+1} \right] S_{01} + \Delta t \left[ (1-\alpha_C)\ \eta_1^{m} + \alpha_C\ \eta_1^{n+1} \right] S_{10}, \\
\eta_1^{n+1} &= \eta_1^{m} - \Delta t \left[ (1-\alpha_C)\ \eta_1^{m} + \alpha_C\ \eta_1^{n+1} \right] S_{10} + \Delta t \left[ (1-\alpha_C)\ \eta_0^{m} + \alpha_C\ \eta_0^{n+1} \right] S_{01}, \\
\end{split}
\end{equation}

where $\alpha_C$ determines whether the transfer terms in the continuity equation are numerically treated explicitly ($\alpha_C=0$) or implicitly ($\alpha_C = 1$). Re-arranging these equations for $\eta_i^{n+1}$, we get

\begin{equation}
\label{transfers_energy3_0}
\begin{split}
\eta_0^{n+1} &= (1-\lambda_{C01}) \eta_0^{m} + \lambda_{C10} \ \eta_1^{m}, \\
\eta_1^{n+1} &= (1-\lambda_{C10}) \eta_1^{m} + \lambda_{C01} \ \eta_0^{m}, \\
\end{split}
\end{equation}

where

\begin{equation}
\begin{split}
\lambda_{Aij} &\equiv \frac{\Delta t S_{ij}}{1 + \alpha_A \Delta t (S_{ji} + S_{ij})} \\
\end{split}
\end{equation}

and $A$ is a label used to identify the transfer coefficients. We use $A=C$ in the continuity equation (mass transfers) and we will later use $A=M$ and $A=T$ for momentum and temperature transfers respectively. The total mass is clearly conserved as $\sum_i \eta_i^{n+1} = \sum_i \eta_i^{n}$ and the total potential energy is also conserved. $\lambda_{ij}$ is between 0 and 1 for all $\alpha_A$ when $\Delta t S_{ij} \leq 1$, meaning $\eta_0$ and $\eta_1$ remain positive. When $\alpha_A=1$, any positive $\Delta t S_{ij} > 0$ may be used. \\

\subsection{Transferring fluid properties - Method 1}
\label{section_transfer_momentum_velocity}
We must also model the transfer of velocity and temperature associated with the re-labelling of mass between fluids described in section \ref{section_transfer_continuity}. The new value of the variable $\phi_i \in [\theta_i, \bm{u}_i]$ should be bounded by the old values of fluids $i$ and $j$ (at time level $m$) so that new extrema are not generated. Also, momentum should be conserved and energy should not increase. Assuming operator-split transfers, the new fluid properties for fluids $0$ and $1$ are written as

\begin{equation}
\begin{split}
\phi_0^{n+1} &= \phi_0^{m} - (1-\alpha_A) \Delta t\ \frac{\eta_1^{q}}{\eta_0^{r}} S_{10} ( \phi_0^{m} - \phi_1^{m} ) - \alpha_A\ \Delta t\ \frac{\eta_1^{q}}{\eta_0^{r}} S_{10} ( \phi_0^{n+1} - \phi_1^{n+1} ), \\
\phi_1^{n+1} &= \phi_1^{m} - (1-\alpha_A) \Delta t\ \frac{\eta_0^{q}}{\eta_1^{r}} S_{01} ( \phi_1^{m} - \phi_0^{m} ) - \alpha_A\ \Delta t\ \frac{\eta_0^{q}}{\eta_1^{r}} S_{01} ( \phi_1^{n+1} - \phi_0^{n+1} ), \\
\end{split}
\end{equation}

where $\phi_i^{m}$ are the values after advection. If $\alpha_A=0$, then $\phi_i$ is treated explicitly and $\alpha_A=1$ means $\phi_i$ is treated implicitly. Note that these equations have additional degrees of freedom in the time-level choice for $\eta_i$, where $q$ and $r$ are the time level choices for the numerator and denominator respectively. $A$ is the label for each governing equation. For the momentum and temperature equations, we will use $A=M$ and $A=T$ respectively. Rearranging for $\phi_i^{n+1}$, we obtain

\begin{equation}
\label{velocity_general_transfer_1}
\begin{split}
\phi_0^{n+1} &= (1-\nu_{A10}^{q,r} )  \phi_0^{m} + \nu_{A10}^{q,r} \phi_1^{m}, \\
\phi_1^{n+1} &= (1-\nu_{A01}^{q,r} )  \phi_1^{m} + \nu_{A01}^{q,r} \phi_0^{m}, \\
\end{split}
\end{equation}

where

\begin{equation}
\label{velocity_transfer_term}
\nu_{Aij}^{q,r} = \frac{ \Delta t S_{ij} \frac{\eta_i^{q}}{\eta_j^{r}} }{1 + \alpha_A \Delta t \left[ S_{ij} \frac{\eta_i^{q}}{\eta_j^{r}} + S_{ji} \frac{\eta_j^{q}}{\eta_i^{r}} \right]}.
\end{equation}

With two degrees of freedom in each of $\alpha_C$, $\alpha_A$, $q$ and $r$, a total of 16 different transfer schemes exist using this method. It is trivial to derive the conservation of momentum/internal energy for the following four schemes:

\begin{enumerate}
	\item $\alpha_C = 0$, $\alpha_M = \alpha_T = 0$ with $q=m$, $r=n+1$.	
	\item $\alpha_C = 0$, $\alpha_M = \alpha_T = 1$ with $q=m$, $r=m$.
	\item $\alpha_C = 1$, $\alpha_M = \alpha_T = 0$ with $q=n+1$, $r=n+1$.
	\item $\alpha_C = 1$, $\alpha_M = \alpha_T = 1$ with $q=n+1$, $r=m$.
\end{enumerate}

\begin{figure}[h!]
	\centering
	
	\begin{tabular}{c l}
		{\large \textbf{a)} $q=m$ and $r=n+1$, including scheme 1} &  \\
		\includegraphics[width=0.55\linewidth]{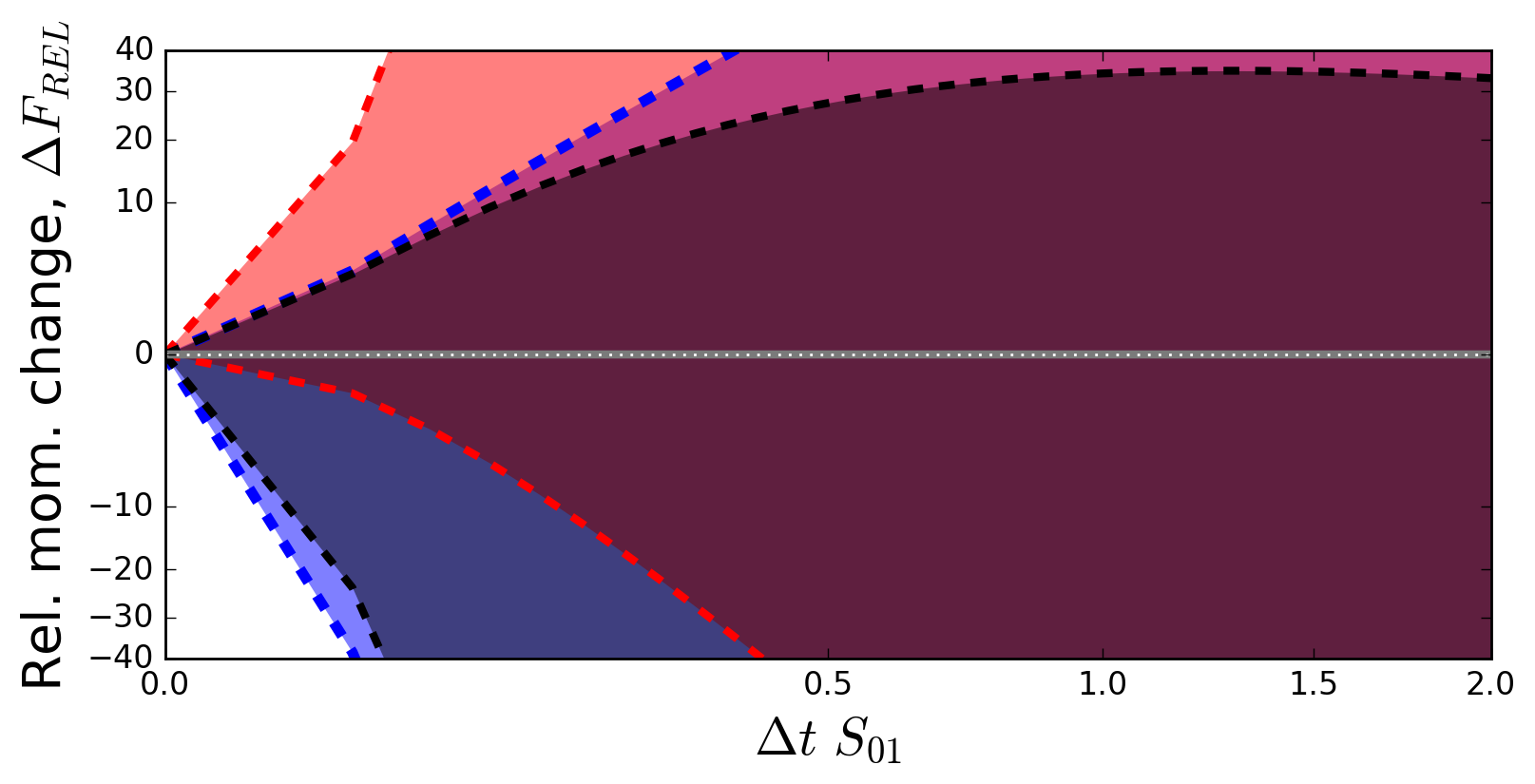} & \includegraphics[width=0.3\linewidth]{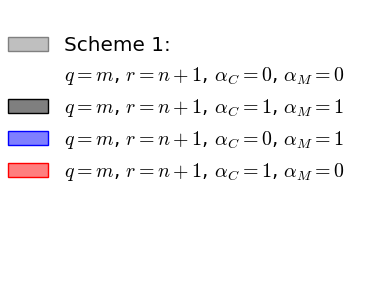} \\
		{\large \textbf{b)} $q=m$ and $r=m$, including scheme 2} & \\ 
		\includegraphics[width=0.55\linewidth]{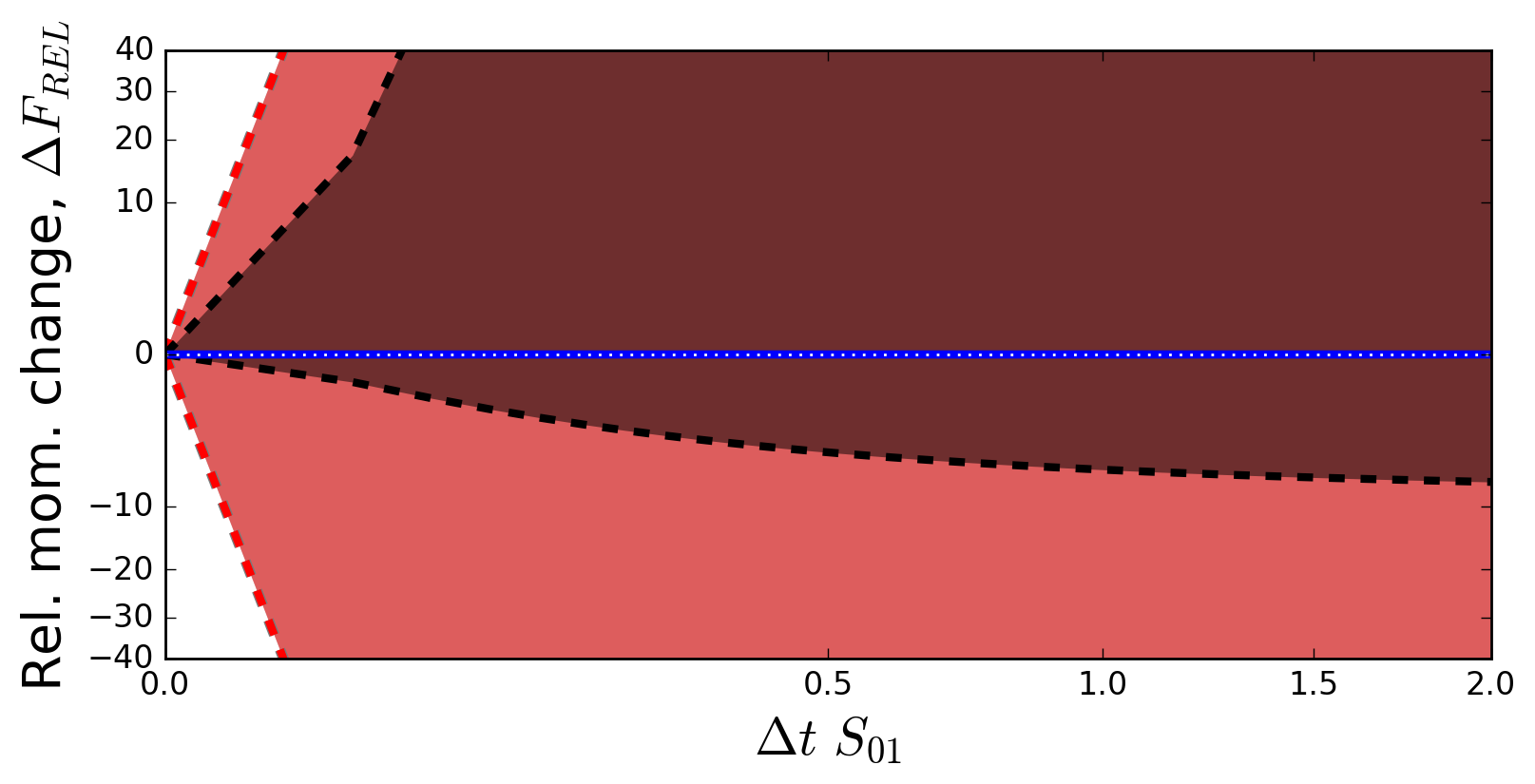} & \includegraphics[width=0.3\linewidth]{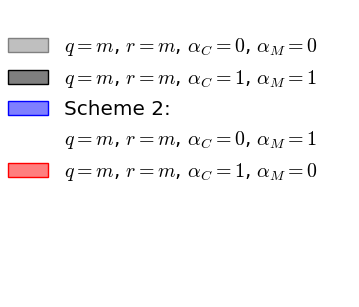} \\
		{\large \textbf{c)} $q=n+1$ and $r=n+1$, including scheme 3} & \\
		\includegraphics[width=0.55\linewidth]{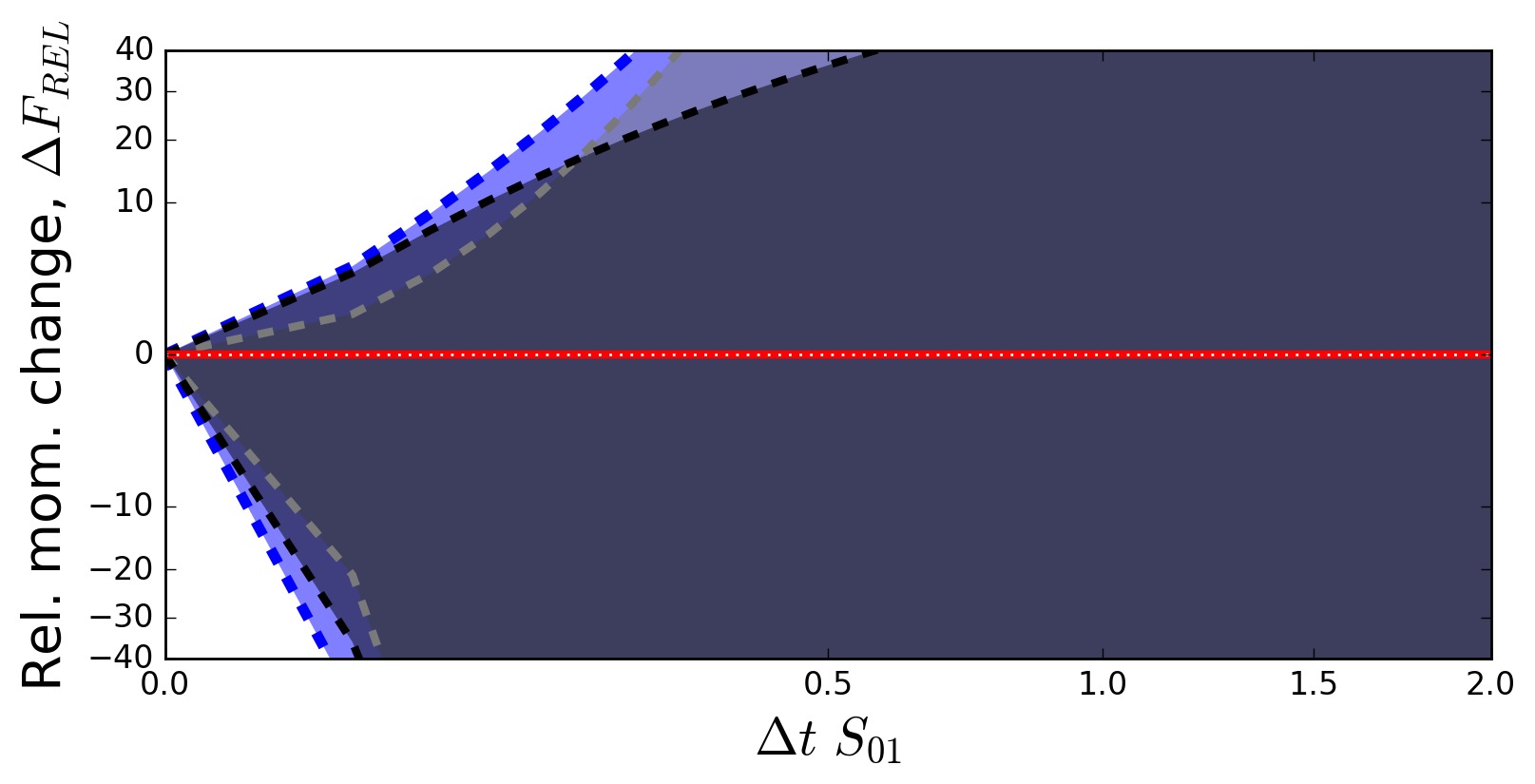} & \includegraphics[width=0.3\linewidth]{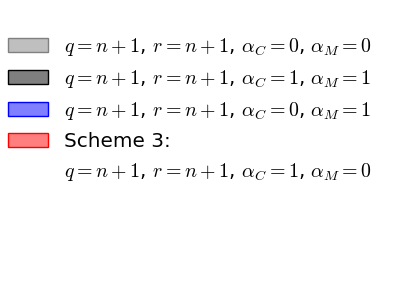} \\
		{\large \textbf{d)} $q=n+1$ and $r=m$, including scheme 4} & \\ 
		\includegraphics[width=0.55\linewidth]{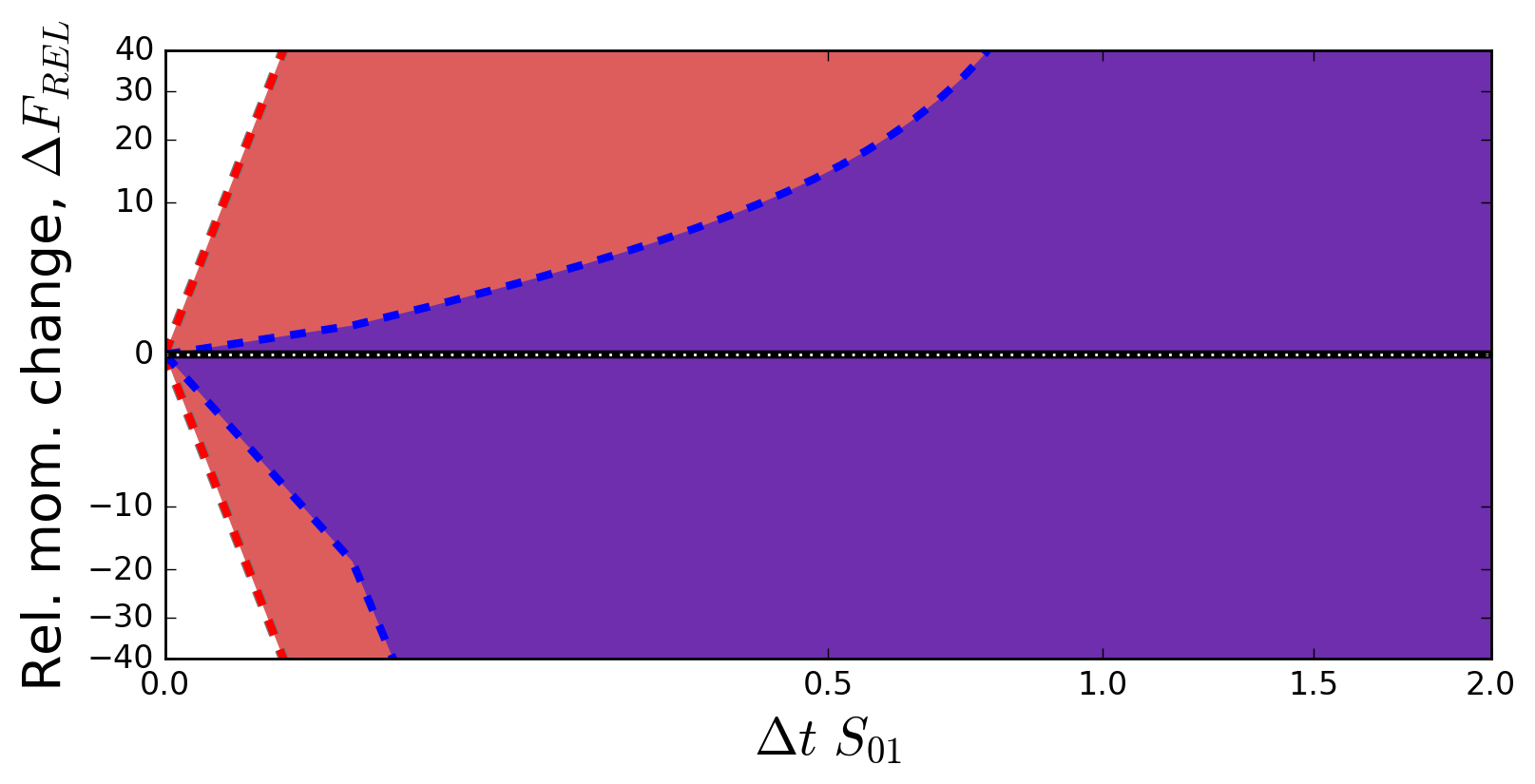} & \includegraphics[width=0.3\linewidth]{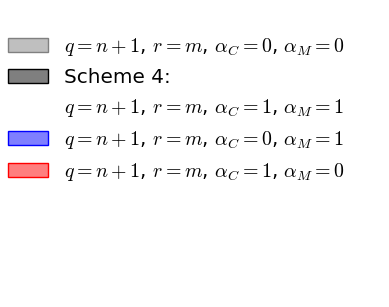} \\
	\end{tabular}
	
	\caption{The momentum changes (relative to initial conditions) of the 16 transfer schemes for method 1. For each scheme, we conduct transfers with initial parameters which include conditions expected in the atmosphere. Schemes are momentum conserving if the relative momentum change due to the transfer is always zero, as indicated by the dotted white line. Schemes which conserve momentum are shown with solid lines whereas dashed lines are used for non-conserving schemes. The conserving schemes are scheme 1 (panel a, grey), scheme 2 (panel b, blue), scheme 3 (panel c, red) and scheme 4 (panel d, black). Square-root scales are used for both axes. }
	\label{figure_momentum_conservation}
\end{figure}

The other 12 schemes do not conserve momentum/internal energy - proof of this is intractable so we instead present numerical analysis of the conservation properties. The relative momentum changes ($\Delta F_\text{REL} \equiv (F^{n+1}-F^{m})/F^{0}$) due to the transfer schemes are calculated using initial conditions which cover a large parameter range, including conditions observed in convective clouds. The transfer schemes were initiated with $\eta_0^m = 1~$kg m$^{-3}$, $u_0^m = 1~$ms$^{-1}$, $\theta_0^m = 300~$K, $\theta_1^m = 301~$K, $S_{01} = 1~$s$^{-1}$. We also use $\Delta t$ in the range $[0,5]~$s, $\eta_1^m$ in the range $[10^{-8},2]~$kg m$^{-3}$, $u_1^m$ in the range $[-150,150]~$ms$^{-1}$ and $S_{10}$ in the range $[0,1]~$s$^{-1}$ - each uniformly discretised 50 times. For a given timestep, $1.25 \times 10^5$ transfers are therefore tested and the range of the relative momentum change for each scheme is plotted. This is shown in figure \ref{figure_momentum_conservation}. These results confirm the momentum conservation analysis of schemes 1-4 (the relative momentum change for these schemes is always zero). The other 12 schemes do not conserve momentum (or internal energy) and will not be analysed further. \\

\begin{figure}[h!]
	\centering
	
	\begin{tabular}{c l}
		{\large \textbf{a)} Method 1 named schemes} &  \\
		\includegraphics[width=0.7\linewidth]{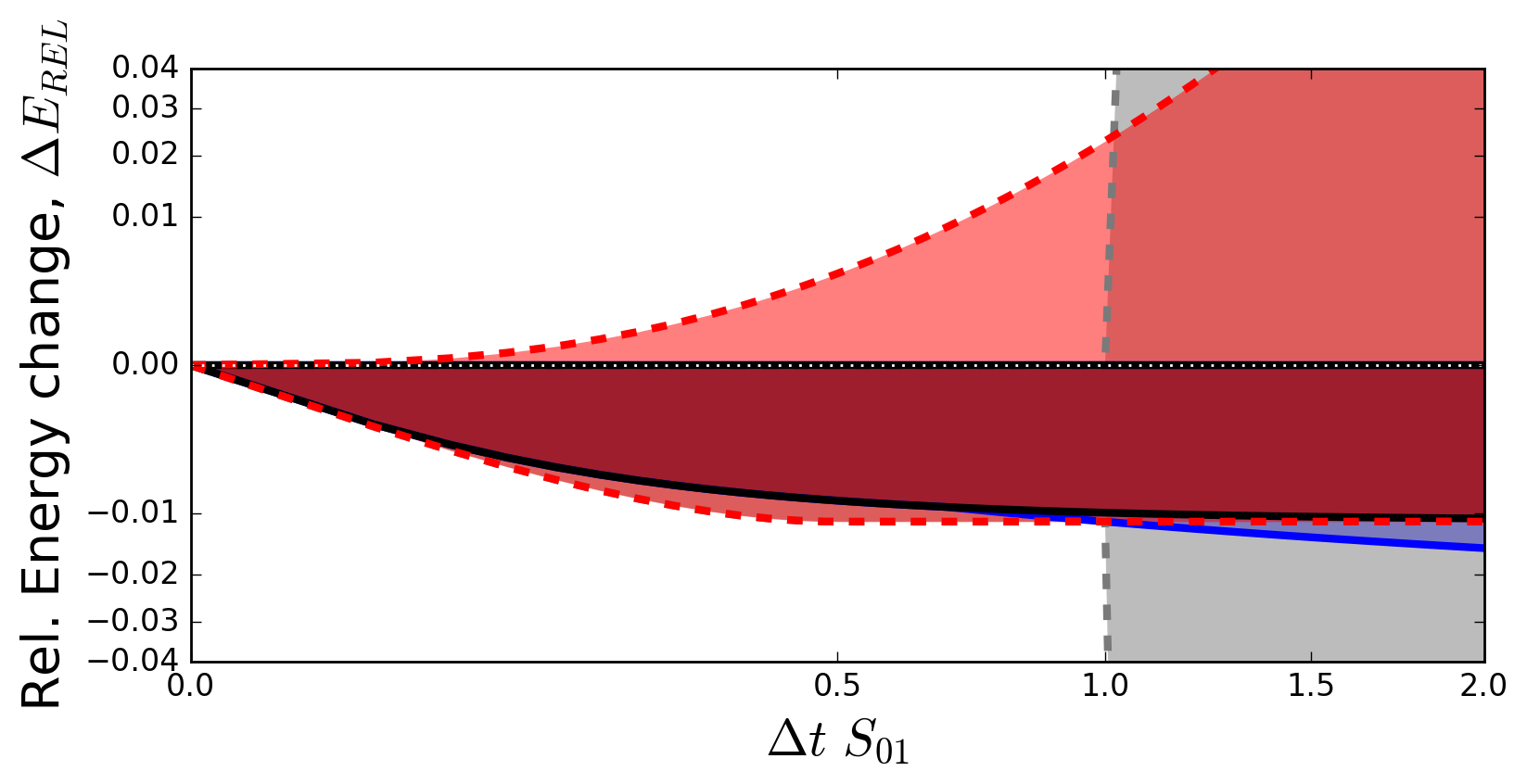} & \includegraphics[width=0.17\linewidth]{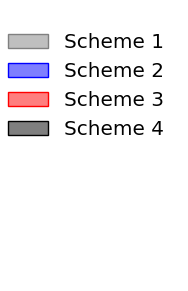} \\
		 & \\
		{\large \textbf{b)} Method 2 named schemes} &  \\
		\includegraphics[width=0.7\linewidth]{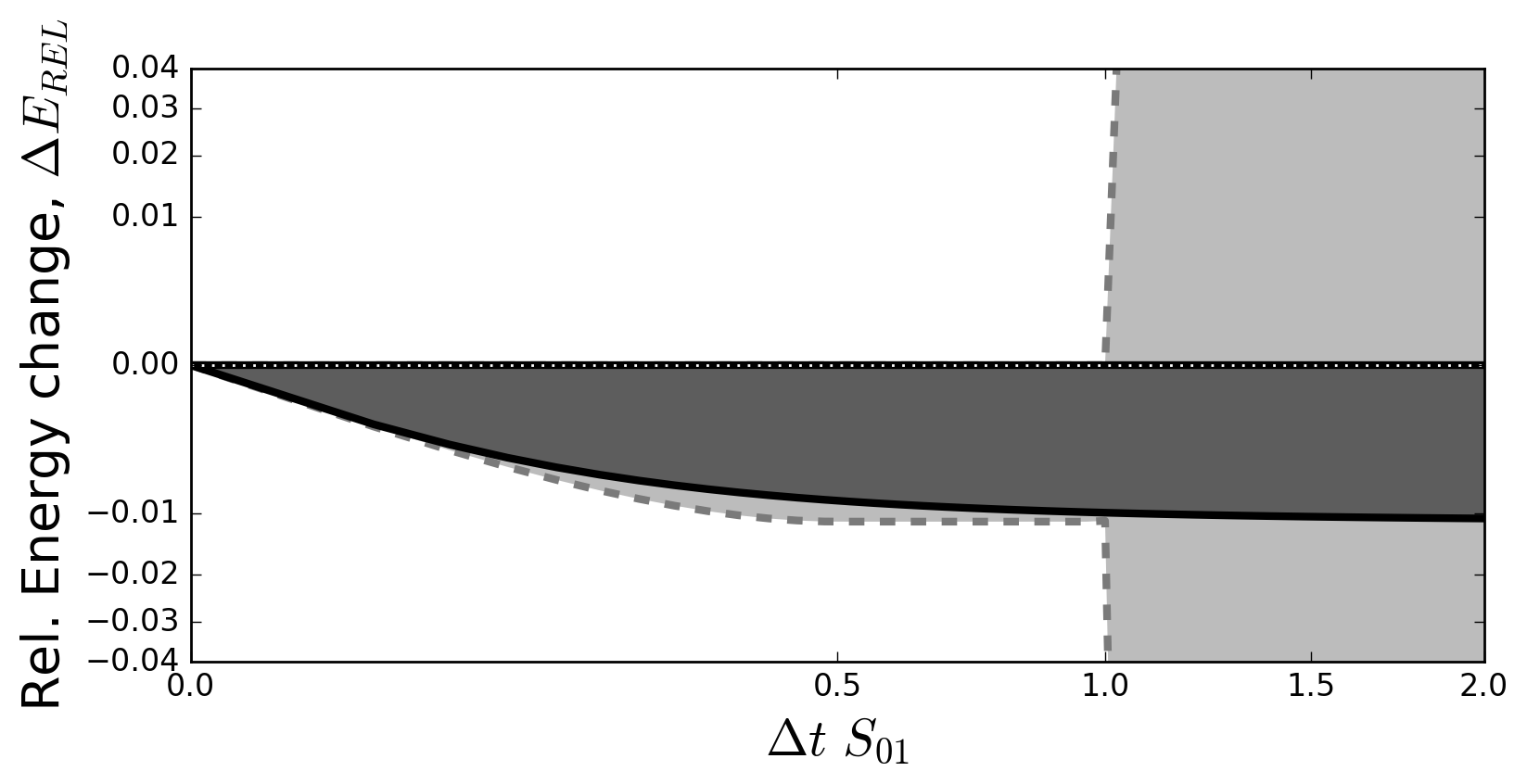} & \includegraphics[width=0.17\linewidth]{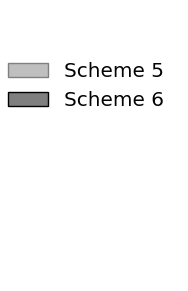} \\
	\end{tabular}
    	
	\caption{The relative energy changes of schemes 1-4 (panel a) and schemes 5 \& 6 (panel b). The minimum and maximum energy changes are calculated using the same parameter-space range as figure \ref{figure_momentum_conservation}. Schemes are energy-diminishing if the relative energy change due to the transfer is never above zero, indicated by the dotted white line. Energy diminishing and energy producing schemes are shown with solid and dashed lines respectively. Scheme 2 (figure a, blue), scheme 4 (figure a, black) and scheme 6 (figure b, black) are energy diminishing. Scheme 1 (figure a, grey) and scheme 5 (figure b, grey) also have energy-diminishing properties for $\Delta t S_{ij} \leq 1$. Square-root scales are used for both axes.}
	\label{figure_energy_conservation}
\end{figure}

For schemes 1-4, $\phi$ will remain bounded if $\nu_{ij}^{q,r} \in [0,1]$, but this can only be guaranteed when $\alpha_A=1$ (we prove this in appendix section \ref{section_transfer_boundedness}) meaning schemes 1 and 3 can produce unbounded velocities and temperatures. Figure \ref{figure_energy_conservation}.a shows the relative energy changes ($\Delta E_\text{REL} \equiv (E^{n+1}-E^{m})/E^{0}$) of schemes 1-4 over the same range of parameter space used in figure \ref{figure_momentum_conservation}. Schemes 2 (blue) and 4 (black) do not increase the total kinetic energy of the system for any $\Delta t S_{ij} > 0$ and scheme 1 (grey) for $\Delta t S_{ij} \leq 1$. Scheme 3 (red) may produce large energy increases for any $\Delta t S_{ij}$. The energy analysis comprehensively samples the parameter-space which is useful for convection modelling, but these results do not concretely prove that schemes 2 and 4 are always energy diminishing. We should therefore also consider other transfer schemes with known energy properties. \\

\subsection{Transferring fluid properties - Method 2 (Mass-weighted transfers)}
\label{section_transfer_momentum_flux}

Transfer terms can also be obtained by considering the flux form equations,
\begin{equation}
\label{momentum_partitioned_flux_transfers}
\frac{\partial (\eta_i \phi_i) }{\partial t} + \bm{\nabla} . ( \eta_i \phi_i \bm{u}_i) =  \eta_i F_i + \sum_{j} \left[ S_{ji} \eta_j \phi_j - S_{ij} \eta_i \phi_i \right],
\end{equation}

which are obtained by combining the continuity equation \ref{continuity_equation} and equation \ref{momentum_partitioned_advective_transfers} with the chain rule. These transfers unconditionally guarantee the conservation of $\eta_i \phi_i$ (momentum and internal energy) as with the mass transfers seen in section \ref{section_transfer_continuity}. By defining our mass-weighted quantity as $\Phi_i \equiv \eta_i \phi_i$, we get:

\begin{equation}
\begin{split}
\Phi_0^{n+1} &= \Phi_0^{m} - \Delta t \left[ (1-\alpha_A)\ \Phi_0^{m} + \alpha_A\ \Phi_0^{n+1} \right] S_{01} + \Delta t \left[ (1-\alpha_A)\ \Phi_1^{m} + \alpha_A\ \Phi_1^{n+1} \right] S_{10}, \\
\Phi_1^{n+1} &= \Phi_1^{m} - \Delta t \left[ (1-\alpha_A)\ \Phi_1^{m} + \alpha_A\ \Phi_1^{n+1} \right] S_{10} + \Delta t \left[ (1-\alpha_A)\ \Phi_0^{m} + \alpha_A\ \Phi_0^{n+1} \right] S_{01}. \\
\end{split}
\end{equation}

The equation takes a similar form to (\ref{mass_transfer_formulation}), meaning we get the solution

\begin{equation}
\label{flux_general_transfer_1}
\begin{split}
\Phi_i^{m} &\equiv \eta_i^{m} \phi_i^{m} \\
\Phi_0^{n+1} &= (1-\lambda_{A01}) \Phi_0^{m} + \lambda_{A10} \ \Phi_1^{m}, \\
\Phi_1^{n+1} &= (1-\lambda_{A10}) \Phi_1^{m} + \lambda_{A01} \ \Phi_0^{m}, \\
\phi_i^{n+1} &= \frac{\Phi_i^{n+1}}{\eta_i^{n+1}}.
\end{split}
\end{equation}

We have proposed this alternative method as we can demonstrate that the total kinetic energy of the system never increases when $\alpha_C = \alpha_A$ (see appendix \ref{section_transfer_energy}). In appendix \ref{section_transfer_boundedness}, we also show that $\phi_i$ is bounded when mass and momentum transfers are treated consistently ($\alpha_C = \alpha_A$) - we will therefore not consider schemes where $\alpha_C \neq \alpha_A$. Using purely explicit or purely implicit treatments, we therefore have 2 more viable transfer schemes for the multi-fluid equations: 

\begin{enumerate}
	\setcounter{enumi}{4}
	\item $\alpha_C = 0$, $\alpha_A = \alpha_M = \alpha_T = 0$.
	\item $\alpha_C = 1$, $\alpha_A = \alpha_M = \alpha_T = 1$.
\end{enumerate}

Note that schemes such as $[\alpha_C = 0.5,\ \alpha_A = 0.5]$ can also be used but there is no increase in order of accuracy as the scheme is operator-split and thus the time level $m$ is not that of the previous timestep. Using time level $n$ instead of $m$ introduces instabilities into the numerical method as updates from the prognostic equations such as advection will be ignored in the transfer scheme. Figure \ref{figure_energy_conservation}.b shows the relative energy changes of schemes 5 and 6 over the same parameter space range used for method 1 schemes. The energy changes are consistent with the analysis in appendix \ref{section_transfer_energy}, whereby scheme 5 never increases in energy for $\Delta t S_{ij} \leq 1$ and scheme 6 for $\Delta t S_{ij} > 0$. \\

\subsection{Transfers on a staggered grid}
\label{section_transfer_experiment_cgrid}
So far, we have assumed that our mass transfers are conducted in the same location, i.e. on a co-located grid (A-grid). But how do the numerical methods change when using staggered grid? Following the C-grid setup used in \cite{weller2018}, we keep our prognostic mass and temperature defined at cell centres and define our velocities on cell faces. Henceforth, a cell-centred variable ($\mu$) which is linearly-interpolated onto cell faces will be denoted by $[\mu]_f$ and a variable defined on cell faces will be denoted by $[\mu]_c$ when it is interpolated onto the cell-centres. \\

The numerical transfer schemes for the mass and potential temperature remain the same, but some adjustments must be made for the velocity transfers (method 1): \\

\begin{equation}
\begin{split}
w_0^{n+1} &= \left( 1-[\nu_{M10}^{q,r}]_f \right)\  w_0^{m} + [\nu_{M10}^{q,r}]_f\ w_1^{m}, \\
w_1^{n+1} &= \left( 1-[\nu_{M01}^{q,r}]_f \right)\  w_1^{m} + [\nu_{M01}^{q,r}]_f\ w_0^{m}, \\
\end{split}
\end{equation}

where $w$ is the vertical velocity. $[\nu_{M10}^{q,r}]_f$ has various degrees of freedom in the choice of interpolations, such as $[S_{ij}]_f \frac{[\eta_i^q]_f}{[\eta_j^r]_f}$ or $\left[ S_{ij} \frac{\eta_i^q}{\eta_j^r} \right]_f$, for example. We will use $\frac{[S_{ij}\ \eta_i^q]_f}{[\eta_j^r]_f}$, once again following \cite{weller2018}. The momentum transfers for method 2 become

\begin{equation}
\begin{split}
N_0^{n+1} &\equiv \left( 1 - [\lambda_{C01}]_f \right) [\eta_{0}^{m}]_f + [\lambda_{C10}]_f \ [\eta_{1}^{m}]_f, \\
N_1^{n+1} &\equiv \left( 1 - [\lambda_{C10}]_f \right) [\eta_{1}^{m}]_f + [\lambda_{C01}]_f \ [\eta_{0}^{m}]_f, \\
\bm{F}_i^{m} &\equiv [\eta_i^{m}]_f \bm{u}_i^{m} \\
\bm{F}_0^{n+1} &= \left( 1 - [\lambda_{M01}]_f \right) \bm{F}_0^{m} + [\lambda_{M10}]_f \ \bm{F}_1^{m}, \\
\bm{F}_1^{n+1} &= \left( 1 - [\lambda_{M10}]_f \right) \bm{F}_1^{m} + [\lambda_{M01}]_f \ \bm{F}_0^{m}. \\
\bm{u}_i^{n+1} &= \frac{\bm{F}_i^{n+1}}{N_i^{n+1}},
\end{split}
\end{equation}

where $[\lambda_{Aij}]_f = \frac{\Delta t [S_{ij}]_f}{1 + \alpha_M \Delta t ([S_{ji}]_f + [S_{ij}]_f)}$. $N_i$ is the fluid mass calculated by conducting the mass transfers on the cell faces - this aids in a consistent and accurate conversion of the mass flux to the fluid velocity. With velocities defined on cell faces, the kinetic energy is calculated on the faces and then interpolated back onto the cell centres:

\begin{equation}
\label{kinetic_energy_cgrid}
E_K^{n+1} = \sum_i \frac{1}{2} \left[ N_i^{n+1} (w_i^{n+1})^2 \right]_c.
\end{equation}

This interpolation method ensures kinetic energy is conserved when converting to the cell centre values \cite[]{ringler2010}.

\subsection{Summary of proposed transfer terms}
We have presented 6 numerical transfer schemes which maintain positivity of mass and conserve mass, momentum, potential energy and internal energy for $\Delta t S_{ij} \leq 1$. These schemes are presented in table \ref{table_transfer_properties}. Schemes 2, 4, 5 and 6 keep the fluid temperatures and velocities bounded, although scheme 5 only does this for $\Delta t S_{ij} \leq 1$. Only schemes 2, 4 and 6 are kinetic-energy-diminishing for all timesteps meaning schemes 1, 3 and 5 can cause numerical instabilities if $\Delta t S_{ij}$ is large. From our analysis, we recommend schemes 4 and 6 as they fulfil all the numerical criteria we have set. Scheme 2 is also a viable scheme if the transfer rate is limited to $\Delta t S_{ij} \leq 1$ to maintain positive mass. Section \ref{section_rising_bubble} will test these schemes on 2D staggered grids. \\

\section{Rising bubble test cases}
\label{section_rising_bubble}

\begin{figure}[h!]
	\centering
	\begin{Large}
	Full bubble test case, potential temperature
	\end{Large}
	\begin{tabular}{c c}
		$\theta$, $t=0~s$ & $\theta$, $t=1000~s$ \\
		\includegraphics[width=0.45\linewidth]{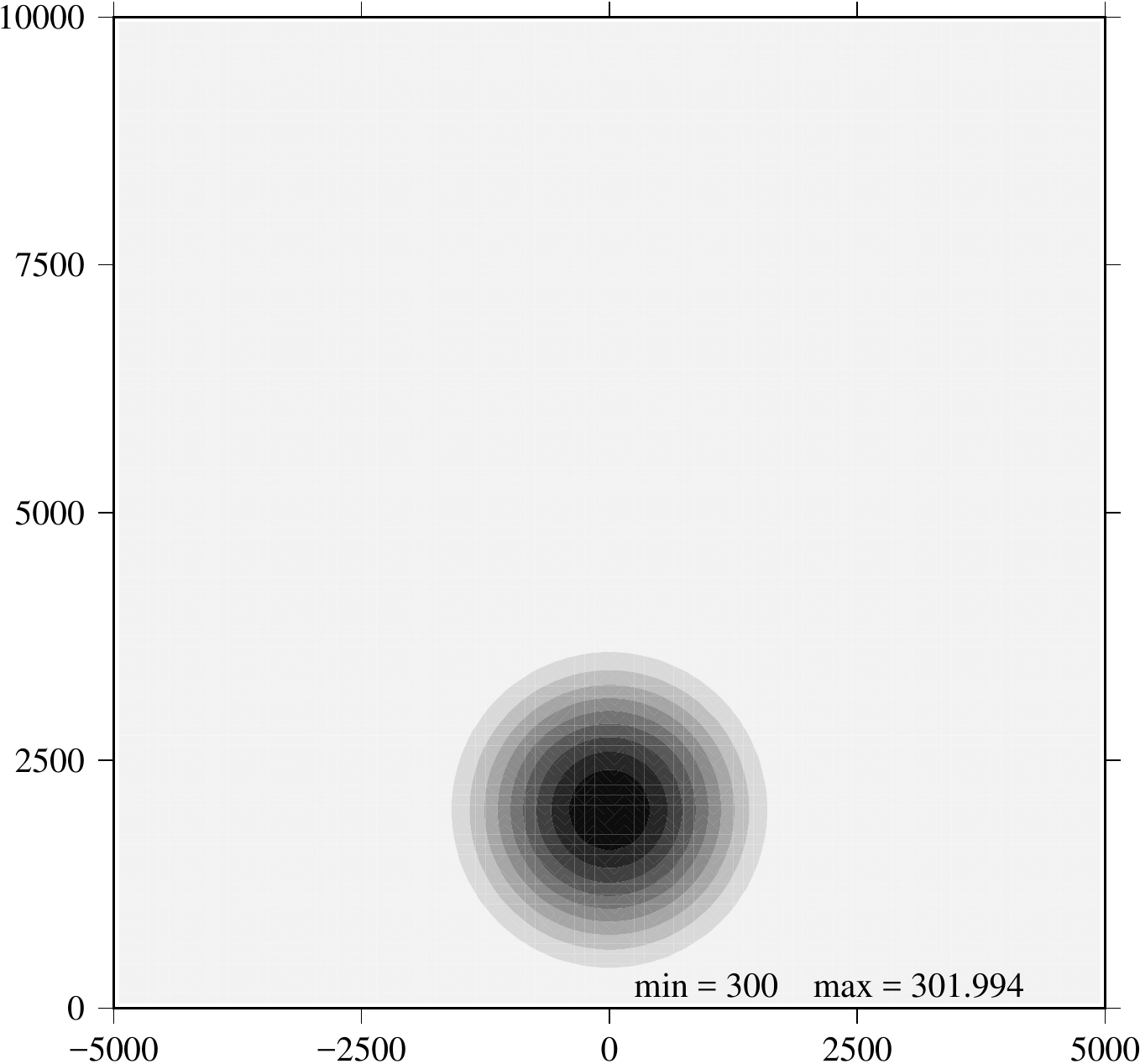} & \includegraphics[width=0.45\linewidth]{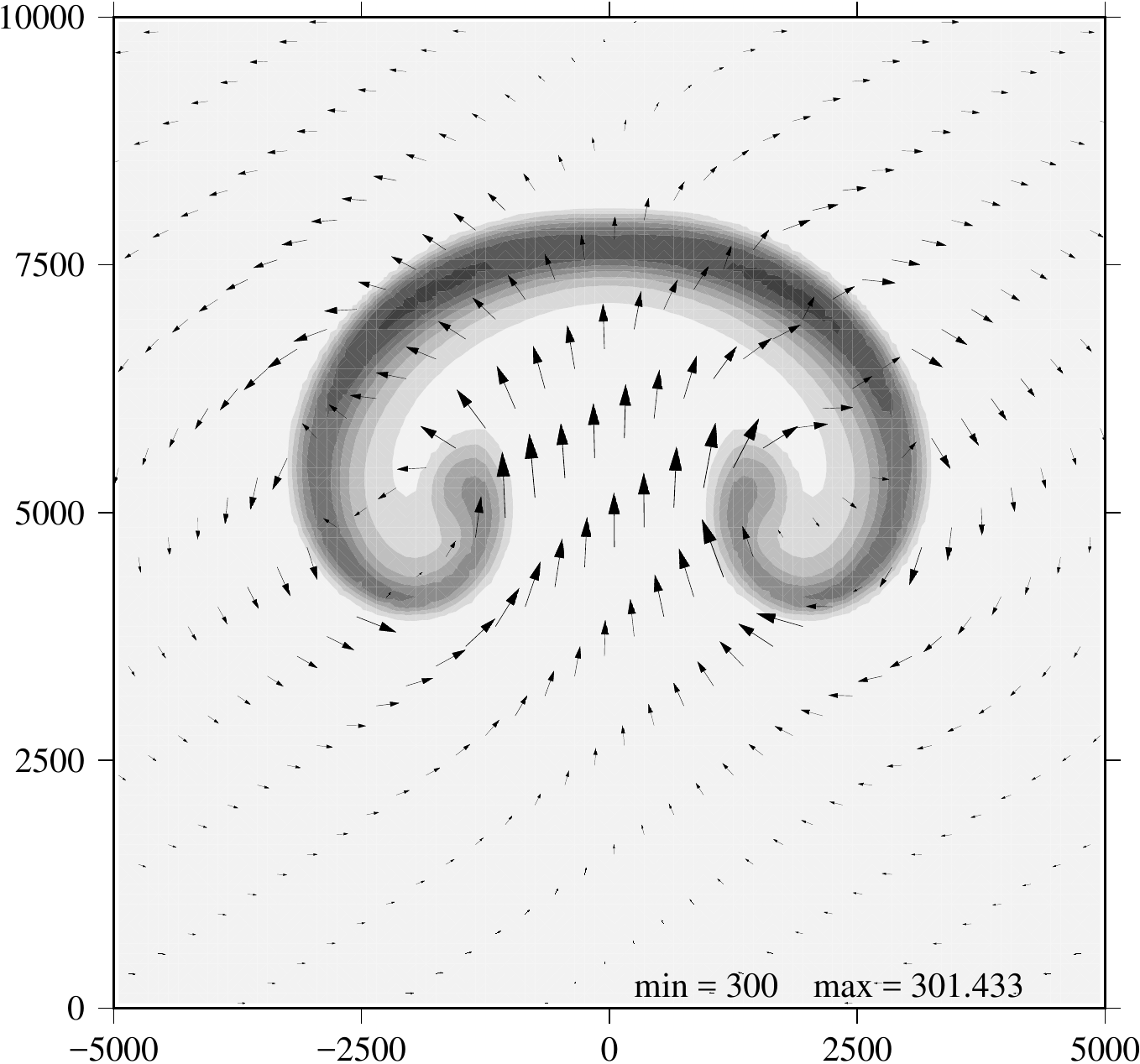}\\
	\end{tabular}
	
	\includegraphics[width=0.95\linewidth]{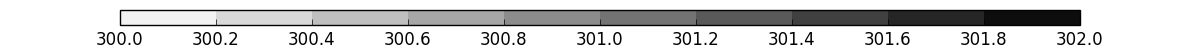} \\
	$\theta (K)$ \\
	\caption{The temperature profile evolution of the full bubble test case, which has the same analytical solution as the single-fluid test case from \cite{bryan2002}. The warm anomaly rises and induces large-scale resolved circulations. The black arrows give the relative magnitudes and directions of the velocity vectors.}
	\label{figure_bubble_transfers}
\end{figure}

In order to test the properties of the various transfer schemes on a staggered grid, we have implemented them into the multi-fluid fully compressible Euler equation solver from \cite{weller2018} using operator splitting. We will run test cases adapted from the single-fluid rising bubble test case \cite[defined in][]{bryan2002} where an initially stationary temperature anomaly rises and generates resolved circulations (see figure \ref{figure_bubble_transfers}). The domain extends to $x \in [-10,10]~$km and $z \in [0,10]~$km with uniform grid spacings $\Delta x = \Delta z = 100~$m and wall boundaries on all sides (where zero-gradient fields are imposed and no fluxes perpendicular to the boundaries). A uniform potential temperature field of $\theta=300~$K is initially chosen with the system in hydrostatic balance and zero velocity. A warm temperature perturbation is then applied at $t=0~$s: \\

\begin{equation}
\theta ' = 2 \cos^2 \left( \frac{\pi}{2} L \right).
\end{equation}

The perturbation is only applied for $L \leq 1$ where $L \equiv \sqrt{\frac{x-x_c}{x_r} + \frac{z-z_c}{z_r}}$, $x_c=10~$km, $z_c=2~$km and $x_r=z_r=2~$km. For the 2-fluid experiments the warm anomaly will be applied to $\theta_1$ only, whereas fluid 0 will remain initialised as $\theta_0=300~$K. \\

We use a 2D C-grid with $\eta_i$ and $\theta_i$ defined at cell centres and the normal component of $\bm{u}_i$ defined at cell faces. The time-stepping is centred Crank-Nicolson with a timestep of $\Delta t = 2~$s.  A van-Leer advection scheme is chosen to maintain positivity of the mass of each fluid. All details of the numerical setup and the numerical solvers used are described in \cite{weller2018}, with the exception of a numerical adjustment which must be made for an operator-split Crank-Nicolson multi-fluid scheme (described in appendix \ref{section_crank_nicolson_adjustments}). \\

\subsection{Full bubble test case}
The first test case is initialised with all mass in fluid 1: $\sigma_0=0$, $\sigma_1=1$. The transfer rate is chosen to transfer a large quantity of fluid 1 to fluid 0:

\begin{align}
S_{01} &= 0, \\
S_{10} &= \frac{1}{\Delta t \ \eta_1^m} \text{max} \left( 0,\ \sigma_\text{min} \rho_0^m - \eta_0^m \right),
\end{align}

where $\sigma_\text{min}=0.1$. This means that the explicit schemes will transfer $10\%$ of the mass in the first timestep (and none thereafter). As fluid 0 initially has no mass, it should inherit the properties of fluid 1 when mass is transferred. We therefore expect the solution to be the same as the single-fluid test case shown in figure \ref{figure_bubble_transfers}. \\

The test case is run for all 20 transfer schemes, including the non-conservative schemes. For each scheme we calculate the relative energy change from the single fluid test case: \\

\begin{equation}
\Delta E_{RSF}^n = \frac{E_{MF}^n - E_{SF}^n}{E_{SF}^0},
\end{equation}

where $E_{SF}^n$ and $E_{MF}^n$ are the total energies at timestep $n$ for the single-fluid and multi-fluid simulations respectively. With a float precision up to 16 decimal places, we expect fluid 0 to inherit the density, velocity and temperature of fluid 1 to machine precision. A relative energy change of $\Delta E_{RSF}^n \sim 10^{-15}$ is therefore reasonable for an energy conserving scheme. \\

The energy changes for all schemes are shown in table \ref{table_onefluid_transfers}. All six conservative schemes produce small energy decreases after one timestep with the implicit mass schemes (3, 4 and 6) producing the smallest energy changes with $\Delta E_{RSF}^1 = -1.18 \times 10^{-15}$ as they have not transferred the full $10\%$ of the mass in the first timestep (unlike schemes 1, 2 and 5). The relative energy changes of schemes 1-6 remain of the order $10^{-15}$ by $t=1000~$s with schemes 5 and 6 having the smallest errors. Many of the non-conservative schemes produce large energy increases due to lack of internal energy conservation and unbounded velocities. Some of these schemes become unstable before the end of the test case at $t=1000~$s. Note that two of the unconservative schemes behave similarly to Schemes 1-6, but internal energy and momentum are not conserved exactly in these schemes. \\

\subsection{Half-bubble test case}

We have already shown solutions for transfers to an empty fluid. But how do the schemes behave when transferring between fluids with comparable mass and different properties? For this we use a 2-fluid test case from \cite{weller2018}, where half the mass is initialised with the warm anomaly (fluid 1) and the other half without (fluid 0):

\begin{equation}
\begin{split}
\sigma_1 &= 
\begin{dcases}
	0.5 & \text{for } L < 1,\\
	0 & \text{otherwise},\\
\end{dcases} \\
\sigma_0 &= 1 - \sigma_1. \\
\end{split}
\end{equation}

The 2-fluid equations with different fluid properties require some form of stabilization \cite[]{weller2018,thuburn2019}. We will use a diffusive mass transfer to couple the fluids:

\begin{equation}
S_{ij} = \frac{1}{2} \frac{K_\sigma}{\eta_i} \text{max} \left( 0,\ \nabla^2 (\eta_j - \eta_i) \right),
\end{equation}

where $K_\sigma = 200~$m$^2~$s$^{-1}$ is a large-enough diffusion coefficient to maintain numerical stability for this test case \cite[as shown in][]{weller2018}. \\

\begin{figure}[h!]
	\centering
	\begin{Large}
	\textbf{a)} Half-bubble test case, potential temperature
	\end{Large}
	\begin{tabular}{c c}	
		$\theta_0$, $t=1000~s$ & $\theta_1$, $t=1000~s$ \\
		\includegraphics[width=0.45\linewidth]{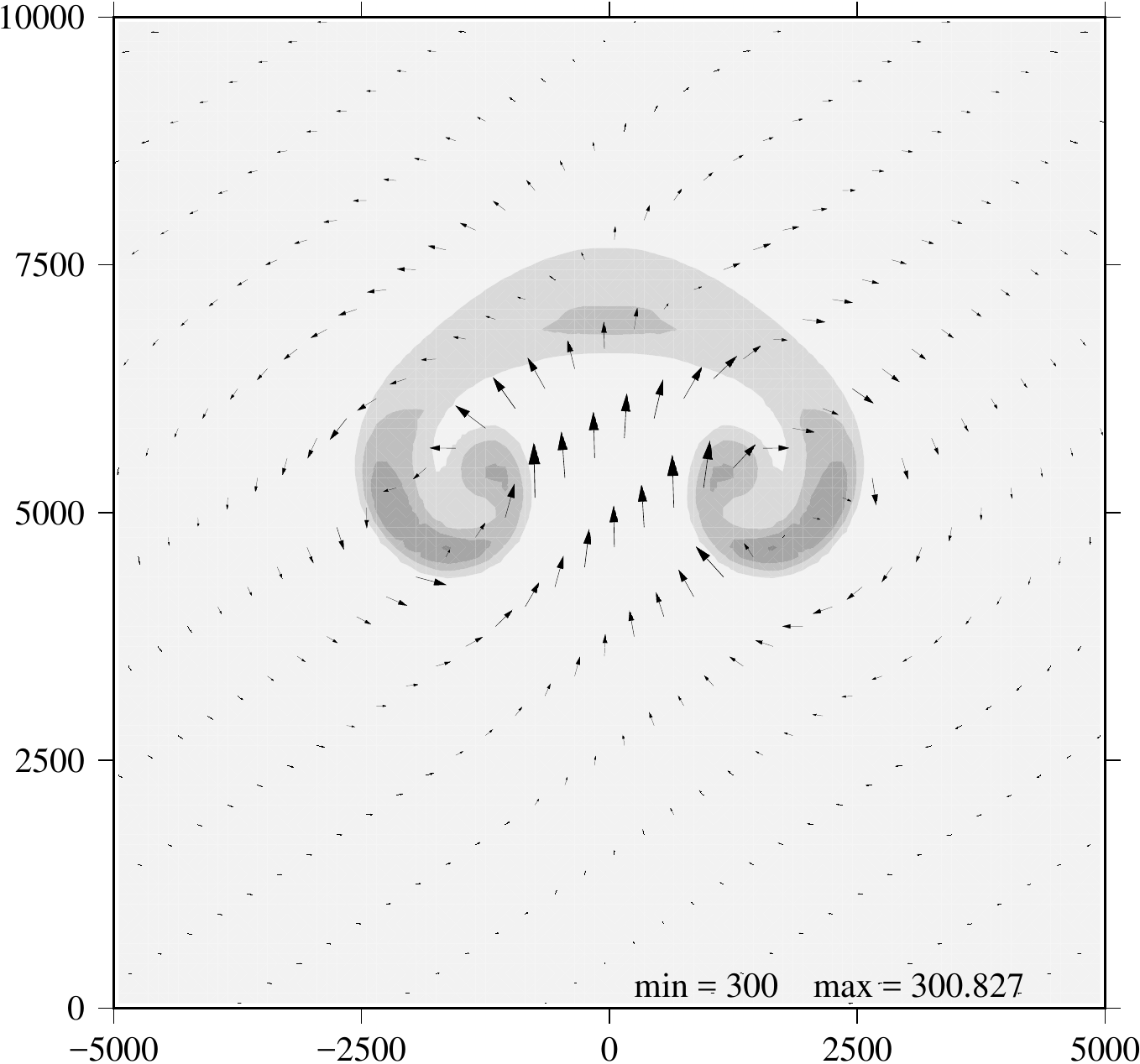} & \includegraphics[width=0.45\linewidth]{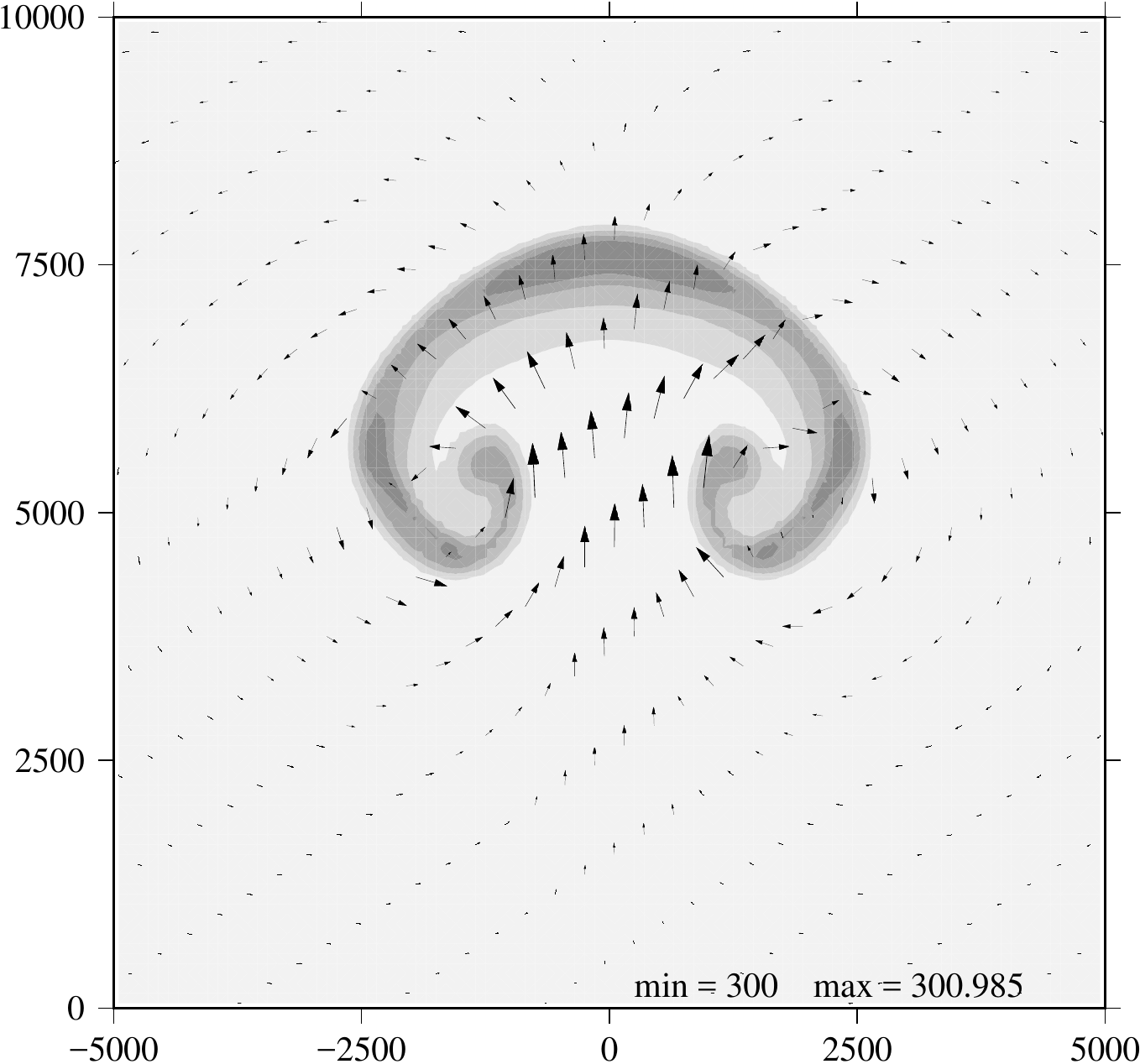}\\
	\end{tabular}
	\includegraphics[width=0.95\linewidth]{bubble_transfers/colorbar_horizontal_theta.png} \\
	$\theta (K)$ \\
	\vspace{10pt}
	\begin{Large}
	\textbf{b)} Half-bubble test case, volume fraction\\
	\end{Large}
	\begin{tabular}{c c}	
		$\sigma_1$, $t=0~s$ & $\sigma_1$, $t=1000~s$ \\
		\includegraphics[width=0.45\linewidth]{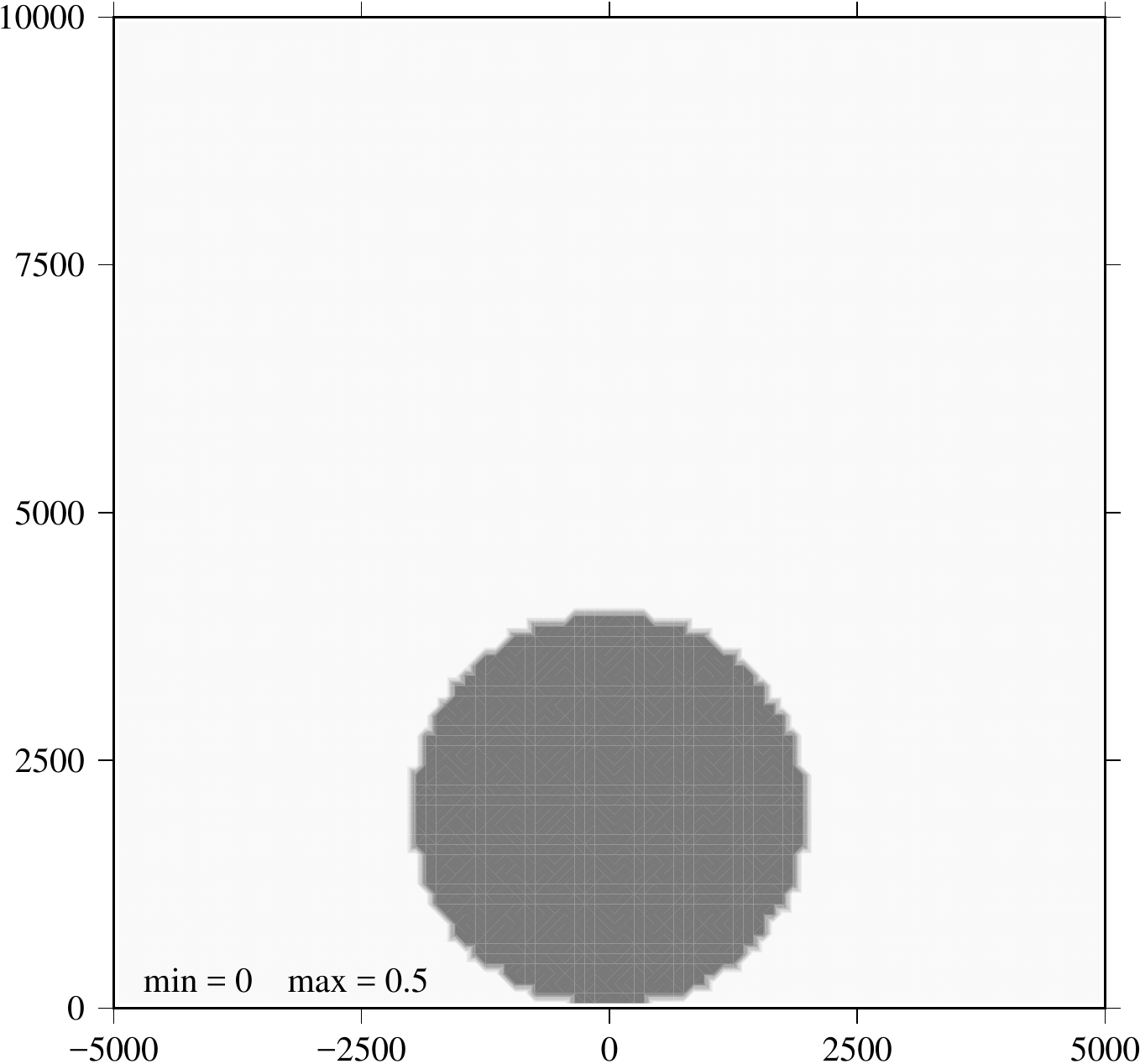} & \includegraphics[width=0.45\linewidth]{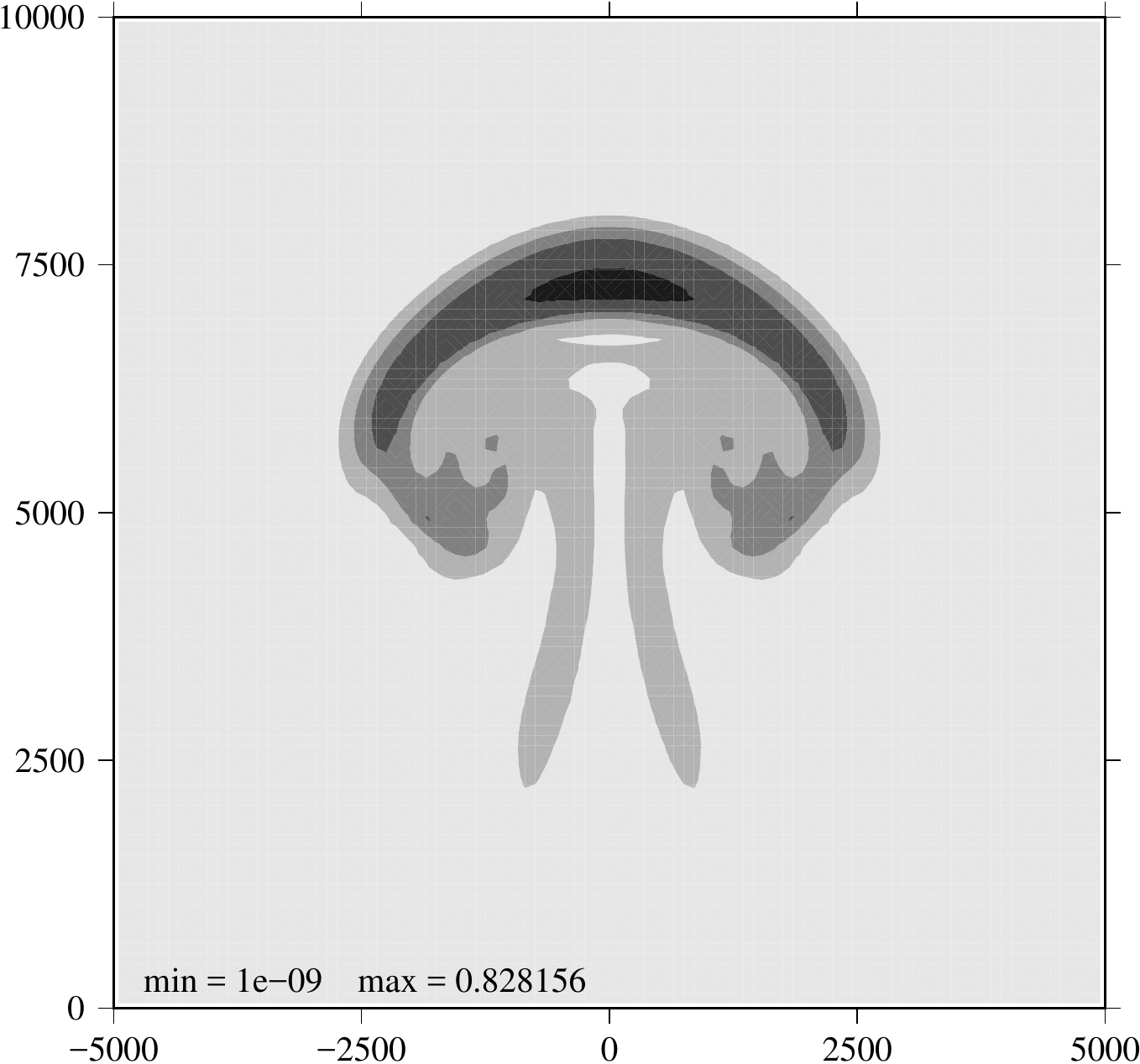}\\
	\end{tabular}	
	\includegraphics[width=0.95\linewidth]{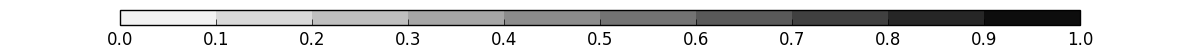} \\
	$\sigma_1$ \\
	\caption{The temperature (a) and volume fraction (b) profiles of the half-bubble test case using scheme 2. The black arrows give the relative magnitudes and directions of the velocity vectors. As the warm anomaly is initially only in half the fluid, the distributions differ from the single-fluid case in figure \ref{figure_bubble_transfers}, including a slower circulation over the domain.}
	\label{figure_bubble_transfers_diff}
\end{figure}

\begin{figure}[h!]
	\centering
	\begin{Large}
	Relative energy change from initial conditions
	\end{Large}	
	\includegraphics[width=\linewidth]{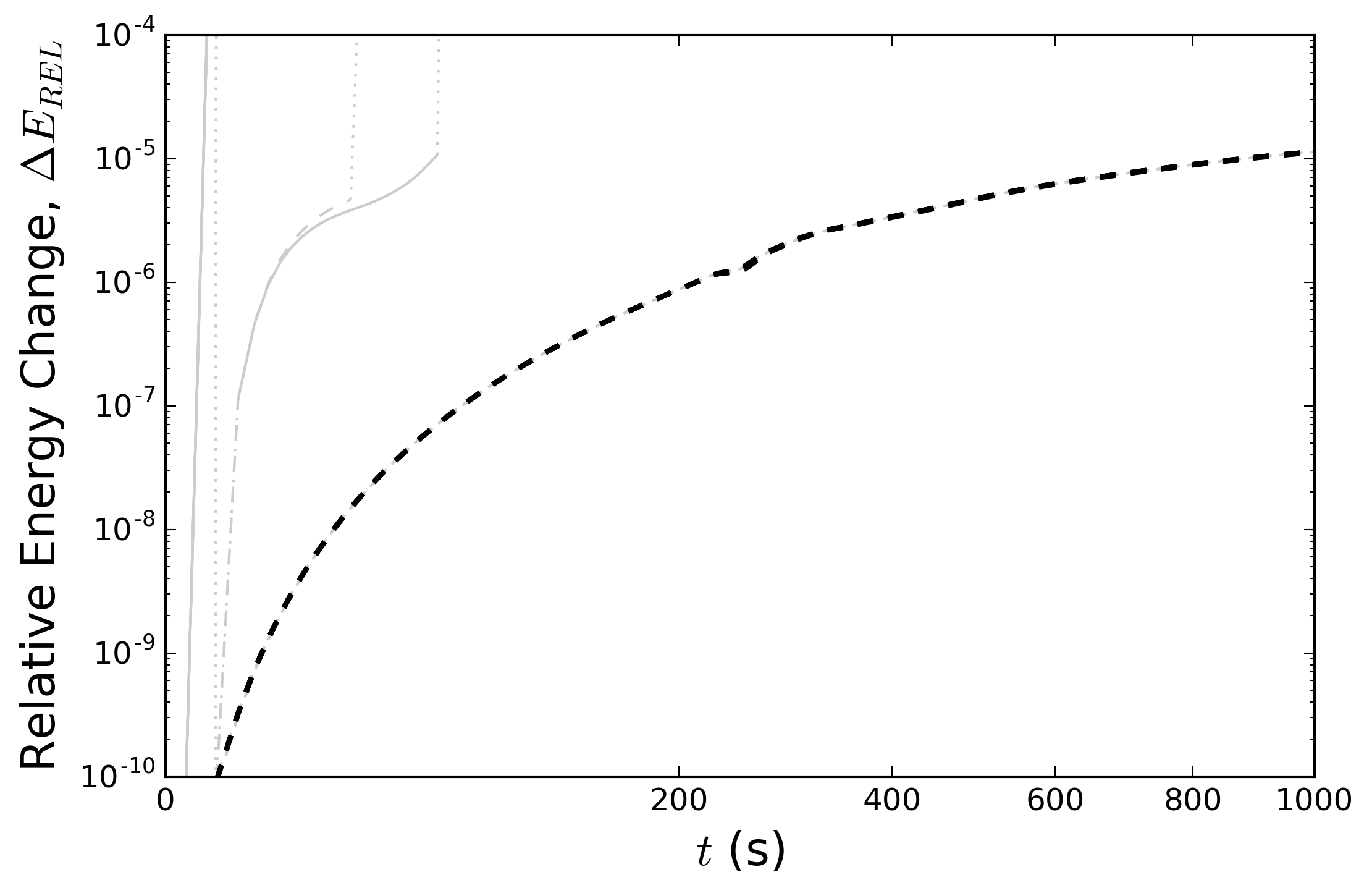}
	\includegraphics[width=0.6\linewidth]{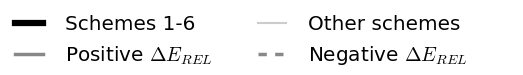}\\
	\caption{The relative energy change from the initial energy for all of the transfer schemes for the 2-fluid rising bubble test case. The solid and dashed lines represent positive and negative energy differences respectively. Dotted lines represent the transition between positive and negative energy changes. Schemes 1-6 are given by the black lines, while the remaining non-conservative schemes are given by the light grey lines.}
	\label{figure_bubble_transfers_2}
\end{figure}

The temperature and volume fraction distributions for this test case are shown in figure \ref{figure_bubble_transfers_diff} - slower circulations form compared with the full bubble test case due to the lower mean temperature anomaly. The energy changes of all schemes relative to the initial conditions are shown in figure \ref{figure_bubble_transfers_2}. Dashed lines represent negative energy changes and solid lines show positive energy changes. Schemes 1-6 follow similar energy evolutions, where energy decreases relative to the initial conditions. The non-conservative schemes (light grey) exhibit various behaviours; many blow up within the first timesteps and produce large energy increases whereas some schemes (which use implicit transfers) follow similar energy evolutions to the conservative schemes. Note that the energy changes due to the numerics of the Crank-Nicolson scheme are far larger than changes due to the transfer schemes - hence all conservative schemes appear to behave similarly. \\

These simulations on a staggered grid are consistent the analysis of the transfer schemes in section \ref{section_transfers} as schemes 1-6 conserve mass, momentum, internal energy and potential energy and are stable for the given test cases with no positive increases in total energy. \\

\section{Conclusions}
Transfer terms between the fluid components in the multi-fluid equations can be used to couple the fluids, represent physical exchanges and stabilise the equations, but the numerical treatment of the transfer terms must be stable. We have presented various numerical methods for treating the transfer terms between the fluids. These schemes are applicable to any multi-fluid equation set where the mean properties of a fluid are transferred with its mass. We have shown that some of the transfer schemes maintain positive mass, keep prognostic variables bounded, conserve momentum, potential \& internal energy and are kinetic energy diminishing. These properties help to keep the overall numerical scheme accurate and stable on co-located and staggered grids. Of the 6 conservative schemes, we have shown that scheme 3 can produce energy increases (figure \ref{figure_energy_conservation}). We have also shown that schemes 5 and 6 do not increase energy for $\Delta t S_{ij} \leq 1$ and $\Delta t S_{ij} > 0$ respectively. We have not proved this for schemes 1, 2 and 4 but have thoroughly explored the relevant parameter space for convection and have not found any instances of kinetic energy increases (other than scheme 1 for $\Delta t S_{ij} > 1$). The fully implicit schemes (schemes 4 and 6) automatically handle large mass transfers and scheme 6 also produces the smallest energy changes in the full bubble test case. As the energy properties are exactly known, scheme 6 has the most desirable numerical properties but scheme 2 (when enforcing mass positivity) and scheme 4 are also good candidates. By using any of these three schemes, the entrainment and detrainment in the multi-fluid equations can be conducted in a numerically stable manner. The physical form of the entrainment and detrainment transfer terms in the multi-fluid equations should be the focus of future studies so that convective processes can be accurately represented. \\

\section{Appendix}

\subsection{Boundedness properties}
\label{section_transfer_boundedness}

For a 2-fluid system, bounded velocity transfer terms can be generalised by

\begin{equation}
\begin{split}
\bm{u}_0^{n+1} &= (1-\beta_{10}) \bm{u}_0^{m} + \beta_{10} \bm{u}_1^{m}, \\
\bm{u}_1^{n+1} &= (1-\beta_{01}) \bm{u}_1^{m} + \beta_{01} \bm{u}_0^{m}, \\
\end{split}
\end{equation}

where $0 \leq \beta_{ij} \leq 1$ ensures the new velocities are bounded. Method 1 has $\beta_{ij} = \nu_{Mij}^{q,r}$ (defined in equation \ref{velocity_transfer_term}). $\nu_{Mij}^{q,r}$ is clearly positive given positive mass and transfer rates. To investigate whether $\nu_{Mij}^{q,r} \leq 1$, we make the denominator small (the worst case scenario) such that $S_{ji}=0~$s$^{-1}$. This gives

\begin{equation}
\begin{split}
\nu_{Mij}^{q,r} = \frac{ \Delta t S_{ij} \frac{\eta_i^{q}}{\eta_j^{r}} }{1 + \alpha_M \Delta t S_{ij} \frac{\eta_i^{q}}{\eta_j^{r}}} &\leq 1 \\
\Delta t S_{ij} \frac{\eta_i^{q}}{\eta_j^{r}} &\leq \frac{1}{1-\alpha_M}. \\
\end{split}
\end{equation}

This boundedness condition is only guaranteed if $\alpha_M = 1$. \\

When $\alpha_M=\alpha_C$, method 2 has $\beta_{ij} = \frac{\lambda_{Cij} \eta_i^{m}}{\lambda_{Cij} \eta_i^{m} + (1-\lambda_{Cji}) \eta_j^{m}}$ which is bounded if $0 \leq \lambda_{Cij} \leq 1$. This is always true for $\Delta t S_{ij} \leq 1$, although boundedness is also guaranteed for $\Delta t S_{ij} > 1$ when $\alpha_C=\alpha_M = 1$.

\subsection{Energy properties}
\label{section_transfer_energy}

For method 2, the total momentum is conserved if the new momenta ($\bm{F}_i^{n+1}$) satisfy 

\begin{equation}
\label{momentum_conservation}
\begin{split}
\bm{F}_0^{n+1} &= (1-\lambda_{M01}) \bm{F}_0^{m} + \lambda_{M10} \bm{F}_1^{m}, \\
\bm{F}_1^{n+1} &= (1-\lambda_{M10}) \bm{F}_1^{m} + \lambda_{M01} \bm{F}_0^{m}, \\
\end{split}
\end{equation}

where $\bm{F}_i^{m} \equiv \eta_i^{m} \bm{u}_i^{m}$. The new kinetic energy after transfers have been applied is

\begin{equation}
\begin{split}
\frac{1}{2} \bm{u}_0^{n+1} . \bm{F}_0^{n+1} + \frac{1}{2} \bm{u}_1^{n+1} . \bm{F}_1^{n+1} &= \frac{1}{2} \bm{u}_0^{m} . \bm{F}_0^{m} + \frac{1}{2} \bm{u}_1^{m} .  \bm{F}_1^{m} - \Delta K, \\
\end{split}
\end{equation}

where $\Delta K \equiv \frac{1}{2} \left( \bm{u}_0^{m} - \bm{u}_1^{m} \right) . \left( \mu_{01} \bm{u}_0^{m} - \mu_{10} \bm{u}_1^{m} \right)$ and $\mu_{ij} \equiv \left[ \lambda_{Mij}(1 - \beta_{ij} - \beta_{ji}) + \beta_{ji} \right] \eta_i^{m}$. When $\mu_{01}=\mu_{10}$ and $\mu_{01} \geq 0$, the kinetic energy will never increase. For method 2 (and when $\alpha_C = \alpha_M$), $\mu_{ij}$ is given by 

\begin{equation}
\mu_{ij} = \frac{\eta_i^{m} \eta_j^{m} \left[ (1-\lambda_{Cij})\ \lambda_{Cij}\ \eta_i^{m} + (1-\lambda_{Cji})\ \lambda_{Cji}\ \eta_j^{m} \right]}{\left[ \lambda_{Cij}\ \eta_i^{m} + (1-\lambda_{Cji})\ \eta_j^{m} \right] \left[ \lambda_{Cji}\ \eta_j^{m} + (1-\lambda_{Cij})\ \eta_i^{m} \right]},
\end{equation}

which is symmetric and always positive meaning this scheme never produces positive energy changes. \\

Such a proof is less trivial for method 1 as these schemes conserve momentum differently compared to equation \ref{momentum_conservation} - $\eta_0 u_1$ and $\eta_1 u_0$ terms are also present with method 1. Instead, the energy changes are calculated over a large range of parameter space and are shown in figure \ref{figure_energy_conservation}.a. \\

\subsection{Numerical adjustments for an operator-split Crank-Nicolson multi-fluid scheme}
\label{section_crank_nicolson_adjustments}
The numerical multi-fluid scheme used for this study follows the implementation by \cite{weller2018}, with the exception of operator-split transfers. For a fluid property such as temperature or velocity ($\phi$) the solution for the Crank-Nicolson scheme (before transfers) is given by

\begin{equation}
\phi_i^{m} = \phi_i^{n} + \Delta t \left[ (1-\alpha) \left( \frac{\partial \phi_i}{\partial t} \right)^{n} + \alpha \left( \frac{\partial \phi_i}{\partial t} \right)^{m} \right],
\end{equation}

where $\alpha$ is the off-centering coefficient. As $\left( \frac{\partial \phi_i}{\partial t} \right)^{n}$ is stored from the previous timestep, we must ensure that it remains consistent with the fluid properties when transfers are made. This is done by computing

\begin{equation}
\begin{split}
\left( \frac{\partial \phi_0}{\partial t} \right)^{n+1} &= \left(1-\nu_{A10}^{q,r} \right)  \left( \frac{\partial \phi_0}{\partial t} \right)^{m} + \nu_{A10}^{q,r} \left( \frac{\partial \phi_1}{\partial t} \right)^{m}, \\
\left( \frac{\partial \phi_1}{\partial t} \right)^{n+1} &= \left(1-\nu_{A01}^{q,r} \right)  \left( \frac{\partial \phi_1}{\partial t} \right)^{m} + \nu_{A01}^{q,r} \left( \frac{\partial \phi_0}{\partial t} \right)^{m}, \\
\end{split}
\end{equation}

for method 1 schemes and 

\begin{equation}
\begin{split}
\left( \frac{\partial \phi_0}{\partial t} \right)^{n+1} &= \frac{(1-\lambda_{A01}) \ \eta_0^{m} \left( \frac{\partial \phi_0}{\partial t} \right)^{m} + \lambda_{A10} \ \eta_1^{m} \left( \frac{\partial \phi_1}{\partial t} \right)^{m}}{\eta_0^{n+1}}, \\
\left( \frac{\partial \phi_1}{\partial t} \right)^{n+1} &= \frac{(1-\lambda_{A10}) \ \eta_1^{m} \left( \frac{\partial \phi_1}{\partial t} \right)^{m} + \lambda_{A01} \ \eta_0^{m} \left( \frac{\partial \phi_0}{\partial t} \right)^{m}}{\eta_1^{n+1}}, \\
\end{split}
\end{equation}

for method 2 schemes. Absence of these terms lead to errors in the numerical solution when using operator-split transfers, especially when a fluid has a small volume fraction or if large transfers are conducted. These terms are not necessary if the Crank-Nicolson off-centering coefficient is set to $\alpha = 1$ but the scheme will be limited to first-order accuracy in time.

\section*{acknowledgements}
The authors acknowledge funding from the NERC RevCon project NE/N013735/1 lead by Bob Plant. RevCon is part of the ParaCon project lead by Alison Stirling at the UK Met Office.

\printendnotes

\begin{table}[h!]
	\begin{center}
	\begin{tabular}{ l l c c c c }
	\hline
	& \textbf{Scheme} & \textbf{Positive $\eta_i$?} & \textbf{Bounded $\theta_i$ \& $\bm{u}_i$?} & \Longunderstack{\textbf{Momentum \&}\\ \textbf{IE conserved?}} & \textbf{KE decreases?} \\
	\hline
	\hline
	& \textbf{\textit{Method 1}} &  &  &  &  \\
	\textbf{1} & $\alpha_C = 0$, $\alpha_A = 0$, $q=m$, $r=n+1$ & \cmark & \xmark  & \cmark\cmark & \cmark \\
	\textbf{2} & $\alpha_C = 0$, $\alpha_A = 1$, $q=m$, $r=m$ & \cmark & \cmark\cmark  & \cmark\cmark & \cmark\cmark \\
	\textbf{3} & $\alpha_C = 1$, $\alpha_A = 0$, $q=n+1$, $r=n+1$ & \cmark\cmark & \xmark  & \cmark\cmark & \xmark \\
	\textbf{4} & $\alpha_C = 1$, $\alpha_A = 1$, $q=n+1$, $r=m$ & \cmark\cmark & \cmark\cmark  & \cmark\cmark & \cmark\cmark \\
	& Other schemes & \cmark & \cmark\ for $\alpha_A=1$  & \xmark & \xmark \\
	& \textbf{\textit{Method 2}} (Mass-weighted transfers) &  &  &  &  \\
	\textbf{5} & $\alpha_C = 0$, $\alpha_A = 0$ & \cmark & \cmark & \cmark\cmark & \cmark \\
	\textbf{6} & $\alpha_C = 1$, $\alpha_A = 1$ & \cmark\cmark & \cmark\cmark  & \cmark\cmark & \cmark\cmark \\
	& Other schemes & \cmark & \cmark\ for $\alpha_C = \alpha_A$  & \cmark\cmark & \cmark\ for $\alpha_C = \alpha_A$ \\
	\hline
	\end{tabular}
	\end{center}
	\caption{The transfer properties of the transfer schemes from methods 1 and 2. Schemes 4 and 6 have all of the ideal transfer properties. In this study, we have not shown that schemes 1 or 4 always decrease energy but we have observed no energy increases in idealised test cases. Ticks indicate that the scheme fulfils the given property for $\Delta t S_{ij} \leq 1$, double ticks show that the property also occurs for all $\Delta t S_{ij} > 0$ and crosses mean the property is not fulfilled.}
	\label{table_transfer_properties}
\end{table}

\begin{table}
	\centering
	
	\begin{tabular}{c c c c c c c c}
	\hline
	\textbf{Name} & \textbf{Method} & $\alpha_C$ & $\alpha_A$ & $q$ & $r$ & $\Delta E_{RSF}^1$ & $\Delta E_{RSF}^{500}$ \\
	\hline
	\hline
	Scheme 1 & 1 & $0$ & $0$ & $m$   & $n+1$ & $\bm{-1.18 \times 10^{-14}}$ & $\bm{-4.71 \times 10^{-15}}$ \\
	$-$      & 1 & $0$ & $1$ & $m$   & $n+1$ & $+6.77 \times 10^{-7}$ & $-1.34 \times 10^{-6}$ \\
	$-$      & 1 & $1$ & $0$ & $m$   & $n+1$ & $-1.23 \times 10^{-7}$ & $-1.97 \times 10^{-6}$ \\
	$-$      & 1 & $1$ & $1$ & $m$   & $n+1$ & $+5.85 \times 10^{-7}$ & $-2.35 \times 10^{-6}$ \\
	\hline
	$-$      & 1 & $0$ & $0$ & $m$   & $m$   & $+1.05 \times 10^{11}$ & $-$ \\
	Scheme 2 & 1 & $0$ & $1$ & $m$   & $m$   & $\bm{-1.18 \times 10^{-14}}$ & $\bm{-7.07 \times 10^{-15}}$ \\
	$-$      & 1 & $1$ & $0$ & $m$   & $m$   & $+9.51 \times 10^{10}$ & $-$ \\
	$-$      & 1 & $1$ & $1$ & $m$   & $m$   & $\bm{-1.18 \times 10^{-15}}$ & $\bm{-6.73 \times 10^{-15}}$ \\
	\hline
	$-$      & 1 & $0$ & $0$ & $n+1$ & $n+1$ & $+1.35 \times 10^{-7}$ & $-6.53 \times 10^{-6}$ \\
	$-$      & 1 & $0$ & $1$ & $n+1$ & $n+1$ & $+7.12 \times 10^{-7}$ & $-1.07 \times 10^{-6}$ \\
	Scheme 3 & 1 & $1$ & $0$ & $n+1$ & $n+1$ & $\bm{-1.18 \times 10^{-15}}$ & $\bm{-6.56 \times 10^{-15}}$ \\
	$-$      & 1 & $1$ & $1$ & $n+1$ & $n+1$ & $+6.15 \times 10^{-7}$ & $-2.00 \times 10^{-6}$ \\
	\hline
	$-$      & 1 & $0$ & $0$ & $n+1$ & $m$   & $+8.47 \times 10^{10}$ & $-$ \\
	$-$      & 1 & $0$ & $1$ & $n+1$ & $m$   & $\bm{-1.18 \times 10^{-14}}$ & $\bm{-4.71 \times 10^{-15}}$\\
	$-$      & 1 & $1$ & $0$ & $n+1$ & $m$   & $+7.86 \times 10^{10}$ & $-$ \\
	Scheme 4 & 1 & $1$ & $1$ & $n+1$ & $m$   & $\bm{-1.18 \times 10^{-15}}$ & $\bm{-8.75 \times 10^{-15}}$ \\
	\hline
	Scheme 5 & 2 & $0$ & $0$ & $-$   & $-$   & $\bm{-1.18 \times 10^{-14}}$ & $\bm{-3.87 \times 10^{-15}}$\\
	$-$      & 2 & $0$ & $1$ & $-$   & $-$   & $+5.13 \times 10^{-14}$ & $-$\\
	$-$      & 2 & $1$ & $0$ & $-$   & $-$   & $+4.58 \times 10^{-14}$ & $-$ \\
	Scheme 6 & 2 & $1$ & $1$ & $-$   & $-$   & $\bm{-1.18 \times 10^{-15}}$ & $\bm{-3.19 \times 10^{-15}}$\\
	\hline
	\end{tabular}
	
	\caption{The relative energy changes of the 2-fluid rising bubble test case, relative to the single-fluid test case. Energy changes are shown for the first timestep ($n=1$) where the transfers are largest and after $1000~$s ($n=500$). Given the machine precision, a relative energy change of the order $10^{-15}$ is expected for energy conservation - results around this range are indicated in bold.}
	\label{table_onefluid_transfers}
\end{table}

\bibliography{bibliography}

\end{document}